\documentclass[smallextended]{svjour3}       
\smartqed  
\usepackage{graphicx}
\usepackage{amsmath, amssymb}
\usepackage{epstopdf, epsfig, multirow, epic, bm}
\usepackage{appendix}  
\usepackage{mathrsfs} 
\usepackage{algorithm,algorithmic} 
\usepackage{array}
\usepackage{makecell}
\usepackage{authblk}
\usepackage{float}
\usepackage{booktabs} 
\usepackage{multirow}
\usepackage{hyperref} 
\hypersetup{
    colorlinks=true, 
    linkcolor=blue,
    urlcolor=blue,  
    citecolor=blue,
}
\usepackage{colortbl}
\usepackage{subfigure}
\usepackage[misc]{ifsym}

\usepackage[table,xcdraw]{xcolor}
\usepackage{listings} 

\makeatletter
\def\@cite#1#2{(#1\if@tempswa , #2\fi)} 
\renewcommand{\@biblabel}[1]{} 
\makeatother

\newtheorem{assumption}{Assumption}

\hypersetup{CJKbookmarks=true}

\journalname{Optimization and Engineering}

\begin{document}
\begin{sloppypar}

\title{An Efficient Quadratic Penalty Method for a Class of Graph Clustering Problems\thanks{The corresponding author's research is supported by NSFC 12071032 and NSFC 12271526}
}

\author{Wenshun Teng \textsuperscript{1} \and
        Qingna Li \textsuperscript{2,~\Letter}
}

\institute{
\begin{itemize}
    \item [\textsuperscript{\Letter}]{Qingna Li (corresponding author)}\\
     \email{qnl@bit.edu.cn}\\
     School of Mathematics and Statistics/Beijing Key Laboratory on MCAACI, Beijing Institute of Technology, Beijing, China
     \at       
     \item[\textsuperscript{1}]{Wenshun Teng}\\
     \email{teng\_wenshun@163.com}\\
     School of Mathematics and Statistics, Beijing Institute of Technology, Beijing, China    
\end{itemize}
}

\date{Received: date / Accepted: date}

\maketitle
\begin{abstract}
Community-based graph clustering is one of the most popular topics in the analysis of complex social networks. This type of clustering involves grouping vertices that are considered to share more connections, whereas vertices in different groups share fewer connections. A successful clustering result forms densely connected induced subgraphs. This paper studies a specific form of graph clustering problems that can be formulated as semi-assignment problems, where the objective function exhibits block properties. We reformulate these problems as sparse-constrained optimization problems and relax them to continuous optimization models. 
We then apply the quadratic penalty method and the quadratic penalty regularized method to the relaxation problem, respectively. Extensive numerical experiments demonstrate that both methods effectively solve graph clustering tasks for both synthetic and real-world network datasets. For small-scale problems, the quadratic penalty regularized method demonstrates greater efficiency, whereas the quadratic penalty method proves more suitable for large-scale cases.
\end{abstract}

\section{Introduction}
\label{intro}
\begin{figure}[htp]
  \centering
  \includegraphics[width=0.5\textwidth]{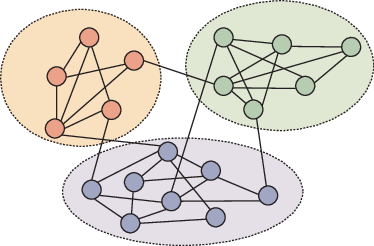}\\
  \caption{A schematic representation of a network with cluster structure. In this network, there are three clusters of densely connected vertices (indicated by solid circles), which are represented by yellow, purple, and green dashed circles. The density of connections between these clusters is much lower}
  \label{fig1}
\end{figure}

The community-based graph clustering problem \cite{no1}, also known as network community detection, is one of the most popular topics in the analysis of complex networks. It has widespread applications in various fields, including sociology \cite{no2,no3}, biology \cite{no4,no5,no6}, and computer science \cite{no7,no8}, with specific uses in areas such as social media \cite{no9}, healthcare \cite{no10}, the web \cite{no11,no12}, and path searches \cite{no13,no14}. A cluster, also referred to as a community, is a set of vertices that are more densely connected to each other than to the rest of the network. The community-based graph clustering problem involves partitioning the vertex set $V$ of a graph into $k$ non-empty subsets, with many edges within each cluster and relatively few edges between clusters. Successful clustering produces clusters with dense internal connections and sparse links between different clusters \cite{no15,GC-BM}. A visual representation of a network with this type of cluster structure is shown in Fig.~\ref{fig1}. It is important to note that this paper focuses on clustering undirected unweighted graphs without self-loops and multiple edges, and we assume that vertices are assigned to non-overlapping clusters.

Unlike the common Euclidean space, graphs do not exist in the form of coordinates, and the distance between two nodes cannot be directly measured. Consequently, traditional clustering algorithms based on Euclidean space, such as K-means, are not applicable to graph clustering problems \cite{no16,GC-BM}. The methods for addressing graph clustering problems can be mainly categorized into five types: (1) spectral algorithms, which utilize the spectral properties of graphs to detect clusters \cite{no17,no18}; (2) statistical inference-based \cite{sta-model} methods, such as Markov \cite{no19} and Bayesian approaches \cite{no20}; (3) dynamics-based methods, which involve running dynamic processes on the network to identify clusters, such as diffusion \cite{no21}, spin dynamics \cite{no22-1,no22} and synchronization \cite{no23}; (4) divisive clustering algorithms, a class of top-down hierarchical methods that recursively partition the graph into clusters \cite{s-1,no24,no25}; (5) optimization-based methods, where the clustering result is obtained by finding the extreme value of a function in the possible clustering space, such as maximizing modularity \cite{no24} (optimizing the modularity quality function) or maximizing similarity \cite{no27}. For a comprehensive and detailed overview of methods for solving graph clustering problems, readers are referred to \cite{s-1,no15,s-3}.

Optimization-based methods are characterized by low computational cost and high accuracy. The maximization of the modularity function \cite{no24} is one of the most popular approaches for community detection, with the Louvain method designed by Blondel et al. \cite{no28} being the most well-known and effective modularity maximization technique. Although these methods offer advantages such as short computation time and the ability to operate without requiring the number of clusters as input parameter \cite{GC-BM}, they also have drawbacks, including the limitations of observed clusters \cite{no29}, the resolution limit problem \cite{no30,cq-2} and the degeneracy problem \cite{no30,cq-1}. Consequently, many researchers have explored alternatives that do not rely on modularity maximization to address graph clustering problems, such as those found in \cite{no33,no34,no35}. Miasnikof et al. \cite{GC-BM} were the first to extend the distance minimization of the binary quadratic formulas of Fan and Pardalos \cite{no27} to the general case of graph clustering problems. 
They employed the Jaccard distance to reflect connectivity and used a Boltzmann machine heuristic to solve the resulting model, which is presented below
\begin{equation}\label{QP}
\begin{aligned}
\min\limits_{x}&\sum\limits_{i}\sum\limits_{j>i}\sum\limits_{k}x_{ik}x_{jk}d_{ij}\\
s.t.&\sum\limits_{k}x_{ik}=1, \ \forall i,\\
&x_{ik}\in \{0, 1\}, \ \forall i, k.
\end{aligned}
\end{equation}
Here $x_{ik}$ is a binary variable that takes the value of $1$ if the vertex $i$ is assigned to the cluster $k$, and $d_{ij}$ represents the Jaccard distance between $i$ and $j$. Note that \eqref{QP} is essentially a binary programming problem with semi-assignment constraints. Inspired by the approach in \cite{no37} to solve hypergraph matching problems, in this paper, we consider a general class of graph clustering problems (see \eqref{semi-assignment}), and propose two methods to solve it.

The contributions of this paper are as follows:
(1) We first summarize several graph clustering models that can be uniformly formulated as semi-assignment problems with objective functions that have block properties. 
(2) For this class of problem, we equivalently reformulate them as sparse-constrained optimization problems. 
(3) We employ both the quadratic penalty method and the quadratic penalty regularized method to solve the relaxation problem. The corresponding subproblems are addressed through an active-set projected Newton method and a spectral projected gradient method, respectively. The special structure of these problems enables computational acceleration in second-order methods.
(4) Extensive numerical results demonstrate that both methods effectively solve graph clustering problems across synthetic graphs and real-world network datasets. For small-scale problems, the quadratic penalty regularized method demonstrates superior computational efficiency, while the quadratic penalty method shows better applicability for complex, large-scale cases.

The remainder of this paper is organized as follows. In Sect.~\ref{sec-2}, we propose the semi-assignment optimization model for graph clustering and discuss its properties. In Sect.~\ref{sec-3}, we investigate the continuous relaxation problem and apply two distinct methods to solve it. In Sect.~\ref{sec-4}, we present extensive numerical results to verify the efficacy of the proposed method. Final conclusions are given in Sect.~\ref{sec-5}.

\section{Semi-assignment optimization model for graph clustering}\label{sec-2}
In this section, we will establish the mathematical model for graph clustering, which can be cased as a semi-assignment optimization problem.
\subsection{Optimization model based on semi-assignment constraints}\label{opt-semi-assignment}
To that end, we consider a graph $G=(V(G), E(G))$, where $V(G)=\{1, 2, \cdots, n\}$ is a nonempty set of vertices and $E(G)=\{(i, j) \ |\ i, j=1, \cdots, n,\ i\ {\rm and }\ j \ {\rm is}\ {\rm connected} \}$ is a set of edges. Let the set of all clusters be $\mathbb{C}= \{C_{1}, C_{2}, \cdots, C_{K}\}$, where $K$ is the number of clusters. We denote the total number of edges (or vertices) by $|E(G)|$ (or $|V(G)|=n$), and the number of vertices in the cluster $i$ by $n^{(i)}$. 

Roughly speaking, graph clustering involves dividing the vertices of a graph into several subsets of densely connected vertices, where vertices within the same subset have more common connections than those in different subsets. Mathematically speaking, let $x_{ik}$ denote whether vertex $i$ is assigned to cluster $k$, and $x_{ik}$ is one if yes and zero otherwise. The essential aim of graph clustering is to assign each vertex a cluster such that some measures of the clustering is minimized. That is,
\begin{equation}\label{semi-assignment}
\begin{aligned}
\min\limits_{x}\hspace{1mm}&f(x)\\
s.t.\hspace{1mm}&\sum\limits_{k=1}^{K}x_{ik}=1, \ i=1,\cdots,n,\\
\hspace{1mm}&x_{ik}\in \{0, 1\}, \ i=1,\cdots,n,\ k=1,\cdots,K,
\end{aligned}
\end{equation}
where $x_{ik}$ is represented in the same way as in \eqref{QP}.
Note that the constraint $\sum\limits_{k}x_{ik}=1$, $x_{ik}\in \{0, 1\}, \ \forall k$, describes the semi-assignment constraint \cite{Integer-programming}. Therefore, \eqref{semi-assignment} is a semi-assignment optimization problem. 

Here, $f(x)$ is a certain criterion to measure the qualification of graph clustering. Various choices of $f$ can be used, such as \cite{no34,GC-BM}. In this paper, we consider the following form of objective function in \cite{GC-BM},
$$f(x)=\sum\limits_{i}\sum\limits_{j>i}\sum\limits_{k}x_{ik}x_{jk}d_{ij},$$
where $d_{ij}$ represents the distance between the vertex $i$ and $j$. Below we mainly consider the following three choices of distances which are Burt's distance \cite{Burt-distance,jaccard-4,jaccard-3}, Jaccard distance \cite{jaccard-1,jaccard-4,jaccard-3} and Otsuka-Ochiai distance \cite{O-O-distance,jaccard-4,jaccard-3}. As will be shown below, this objective function exhibits block properties.

To conclude this subsection, we would like to mention that there are also other criteria for $f$. For example, \cite{no34} focuses on maximizing the mean intra-cluster density while introducing a penalty function $P_{k}(M)$ to discourage clusters that are too large or too small. The penalty functions can be defined as $P_{k}(M)=\max\{0,\sum\limits_{i}x_{ik}\}$ or $P_{k}(M)=(\sum\limits_{i}x_{ik}-M)^{2}$, where $M$ is a parameter. Then, the objective function becomes:
$$f(x)=-\sum\limits_{k=1}^{|\mathbb{C}|}\left[\sum\limits_{i,j}\left(\frac{w_{i,j}x_{ik}x_{jk}}{0.5\times n^{(k)}(n^{(k)}-1)}-\lambda P_{k}(M)\right)\right],$$
where $w_{i,j}$ is the weight of the edge connecting vertex $i$ and vertex $j$, and $\lambda$ is a penalty coefficient.
Unfortunately, this objective function does not exhibit block properties (Prop.~\ref{prop2}) and is computationally challenging. 

\subsection{Different distance metrics}\label{distance}
As mentioned above, graph clustering is defined as subsets of vertices that are considered similar to some extent. This similarity is manifested through the number of shared connections and is translated into a distance metric. Therefore, the distance measure we need is based on similarity rather than the shortest path distance. 
The following lists three types of distance that are often used as similarity measures in graph clustering problems.

\paragraph{Burt's distance $D^{B}$}
The Burt distance \cite{jaccard-4,jaccard-3}, borrowed from sociology \cite{Burt-distance}, is defined as follows for the distance between the vertex $i$ and $j$: 
$$D^{B}_{ij}=\sqrt{\sum\limits_{s\neq i,j}(A_{is}-A_{js})^{2}},$$
where $A=(A_{ij})\in\mathbb{R}^{n\times n}$ is the adjacency matrix. For an unweighted graph, $A_{ij}=1$ if $(i,j)\in E(G)$ and $0$ otherwise. For a weighted graph, if $(i,j)\in E(G)$, $A_{ij}$ is the weight of the edge between $i$ and $j$; otherwise, $A_{ij}=0$. Then Burt's distance matrix is defined as $D^{B}=(D^{B}_{ij})\in\mathbb{R}^{n\times n}$.

\paragraph{Jaccard distance $D^{J}$}
The Jaccard distance \cite{jaccard-4,jaccard-3}, originating from botany \cite{jaccard-1}, is defined between the vertex $i$ and $j$ as follows:
$$D^{J}_{ij}=\left\{\begin{array}{l}
		1-\frac{|a_{i}\cap a_{j}|}{|a_{i}\cup a_{j}|},\ i\neq j,\\
		0,\hspace{1.3cm} i= j,
		\end{array}	\right.$$
where $a_{i}$ denotes the set of all vertices that share an edge with vertex $i$. The ratio $\frac{|a_{i}\cap a_{j}|}{|a_{i}\cup a_{j}|}$ represents the Jaccard similarity. The above applies to the case of an unweighted graph. For a weighted graph, we use the following expression:
$$D^{J}_{ij}=\left\{\begin{array}{l}
		1-\frac{\sum\limits_{s=1}^{n}\min\{w_{is},w_{js}\}}{\sum\limits_{s=1}^{n}\max\{w_{is},w_{js}\}},\ i\neq j,\\
		0,\ i= j,
		\end{array}	\right.$$
where $w_{is}$ denotes the weight of the edge between the vertex $i$ and the vertex $s$. Note that $D^{J}_{ij}\in[0,1]$. Then Jaccard distance matrix is defined as $D^{J}=(D^{J}_{ij})\in\mathbb{R}^{n\times n}$.

\paragraph{Otsuka-Ochiai distance $D^{O}$}
The Otsuka-Ochiai distance originates from zoology \cite{O-O-distance}. 
In the case of unweighted gaphs, Otsuka-Ochiai distance \cite{jaccard-4,jaccard-3} between the vertex $i$ and $j$ is defined as follows: 
$$D^{O}_{ij}=1-\frac{|a_{i}\cap a_{j}|}{\sqrt{|a_{i}|\times |a_{j}|}}\in[0,1],$$
where $a_{i}$ is the same as that in the Jaccard distance. In the weighted case, we use the following expression:
$$D^{O}_{ij}=1-\frac{\sum\limits_{s=1}^{n}\min\{w_{is},w_{js}\}}{\sqrt{\sum\limits_{s=1}^{n}w_{is}\times \sum\limits_{s=1}^{n}w_{js}}}\in[0,1],$$
where $w_{is}$ denotes the weight of the edge between the vertex $i$ and $s$. Then Otsuka-Ochiai distance matrix is defined as $D^{O}=(D^{O}_{ij})\in\mathbb{R}^{n\times n}$.

The three $n\times n$ distance matrices (denoted as $D^{B}$, $D^{J}$ and $D^{O}$ respectively) satisfy Prop.~\ref{prop-jaccard}. 

\begin{proposition}\label{prop-jaccard}
For the matrix $D$ defined by $D^{B}$, $D^{J}$ or $D^{O}$ as above, it holds that $D_{ij}\geq 0$, $D_{ii}= 0$, $D_{ij}=D_{ji}$, $i=1,\cdots,n$, $j=1,\cdots,n$.
In other words, $D$ is nonnegative, symmetric and has zero diagonal elements.
\end{proposition}

\subsection{Equivalent formulation of \texorpdfstring{$\eqref{semi-assignment}$}{}{} }\label{equivalent}
Based on Sect.~\ref{opt-semi-assignment} and Sect.~\ref{distance}, the graph clustering model we consider in this paper takes the following form
\begin{equation}\label{QP1}
\begin{aligned}
\min\limits_{x}&\sum\limits_{i}\sum\limits_{j>i}\sum\limits_{k}x_{ik}x_{jk}d_{ij}\\
s.t.&\sum\limits_{k}x_{ik}=1, \ i=1,\cdots,n,\\
&x_{ik}\in \{0, 1\}, \ i=1,\cdots,n,\ k=1,\cdots,K.
\end{aligned}
\end{equation}

Denote
${\bf x}=({\bf x}^{\top}_{1}, {\bf x}^{\top}_{2},\cdots, {\bf x}^{\top}_{n})^{\top}\in \mathbb{R}^{nK}$, ${\bf x}^{\top}_{i}=(x_{i1}, x_{i2}, \cdots, x_{iK})\in\mathbb{R}^{K}$. 
Then the sum of distances generated by the vertex $i$ and $j$ can be expressed as the following equation
\begin{equation}\label{reform-1}
\sum\limits_{k}x_{ik}x_{jk}d_{ij}=
\left[\begin{array}{c c c c}
x_{i1} \ x_{i2} \ \cdots \ x_{iK}
\end{array}\right]\left[\begin{array}{c c c c}
d_{ij}&0&\cdots&0\\
0&d_{ij}&\cdots&0\\
\vdots&\vdots&\ddots&\vdots\\
0&0&\cdots&d_{ij}
\end{array}\right]
\left[\begin{array}{c}
x_{j1}\\
x_{j2}\\
\vdots\\
x_{jK}
\end{array}\right]=
{\bf x}_{i}^{\top}(d_{ij}I_{k}){\bf x}_{j},
\end{equation}
where $I_{k}\in\mathbb{R}^{K\times K}$ is an identity matrix. Let $\overline{\Delta}_{ij}=d_{ij}I_{k}\in\mathbb{R}^{K\times K}$ and $\Delta=D\otimes I_{k}\in\mathbb{R}^{nK\times nK}$. Then the objective function of \eqref{QP1} can be written as
\begin{equation}\label{reform-2}
\sum\limits_{i}\sum\limits_{j>i}\sum\limits_{k}x_{ik}x_{jk}d_{ij}=\frac{1}{2}
\left[\begin{array}{c c c c}
{\bf x}_{1}^{\top} \ {\bf x}_{2}^{\top} \ \cdots \ {\bf x}_{n}^{\top}
\end{array}\right]
\left(
D\otimes I_{k}\right)
\left[\begin{array}{c}
{\bf x}_{1}\\
{\bf x}_{2}\\
\vdots\\
{\bf x}_{n}
\end{array}\right]=
\frac{1}{2}{\bf x}^{\top}\Delta {\bf x}.
\end{equation}

Hence, \eqref{QP1} can be written in the following equivalent form
\begin{equation}\label{QP2}
\begin{aligned}
\min\limits_{{\bf x}}& \hspace{2mm}f({\bf x}):=\frac{1}{2}{\bf x}^{\top}\Delta {\bf x}\\
{\rm s.t.}&\hspace{2mm}{\bf x}_{i}^{\top}{\bf 1}_{k}=1,\ i=1,\cdots,n,\\
&\hspace{2mm}\textbf{x}\in \{0, 1\}^{nK},
\end{aligned}
\end{equation}
where ${\bf 1}_{k}=[1, 1, \cdots, 1]^{\top}\in\mathbb{R}^{K}$.
\eqref{QP2} can also be written in the following equivalent form:
\begin{equation}\label{QP6}
\begin{aligned}
\min\limits_{{\bf x}\in\mathbb{R}^{nK}}& \hspace{2mm}f({\bf x})\\
{\rm s.t.}&\hspace{2mm}{\bf x}_{i}^{\top}{\bf 1}_{k}=1,\ i=1,\cdots,n,\\
&\hspace{2mm}{\bf x}\geq0,\\
&\hspace{2mm}\|{\bf x}\|_{0}\leq n.
\end{aligned}
\end{equation}
This is because from ${\bf x}_{i}^{\top}{\bf 1}_{k}=1,\ i=1,\cdots,n$ and ${\bf x}\geq0$ we know that $\|{\bf x}\|_{0}\geq n$, and since we have $\|{\bf x}\|_{0}\leq n$, it follows that $\|{\bf x}\|_{0}=n$. Therefore, $\textbf{x}\in \{0, 1\}^{nK}$.

In \eqref{QP6}, we remove the last constraint ${\bf x}\in \{0, 1\}^{nK}$ of \eqref{QP2}, and add two constraints, $\|{\bf x}\|_{0}\leq n$ and ${\bf x}\geq0$.
Due to the constraint ${\bf x}_{i}^{\top}{\bf 1}_{k}=1$, $i=1,\cdots,n$, we know that $\|{\bf x}\|_{0}\geq n$. 
Additionally, because of the last constraint $\|{\bf x}\|_{0}\leq n$, we have $\|{\bf x}\|_{0}= n$. 

The following result is obvious due to the definition of $\Delta$.
\begin{proposition}\label{prop1}
Recall $\Delta=D\otimes I_{k}$. It holds that $\Delta_{ii}=0$, $\Delta_{ij}=\Delta_{ji}$, $i=1,2,\cdots, nK$, $j=1,2,\cdots,nK$.
In other words, $\Delta$ is nonnegative, symmetric and has zero diagonal elements.
\end{proposition}

To present the property of $f(x)$, let an index $i_{0}\in\{1, 2, \cdots, n\}$ and the set $I^{-i_{0}}=\{1, 2, \cdots, i_{0}-1, i_{0}+1, \cdots, n\}$. We can rewrite $f({\bf x})$ in \eqref{QP2} as:
\begin{equation}\label{QP4}
\begin{aligned}
f({\bf x})&=\frac{1}{2}\sum\limits_{i=1}^{n}\sum\limits_{j\neq i}{\bf x}_{i}^{\top}\overline{\Delta}_{ij}{\bf x}_{j}\\
&=\frac{1}{2}\left(2{\bf x}_{i_{0}}^{\top}\sum\limits_{j\neq i_{0}}\overline{\Delta}_{i_{0}j}{\bf x}_{j}+\sum\limits_{i\in I^{-i_{0}}}\sum\limits_{j\neq i, j\in I^{-i_{0}}}{\bf x}_{i}^{\top}\overline{\Delta}_{ij}{\bf x}_{j}\right)\\
&={\bf x}_{i_{0}}^{\top}\sum\limits_{j\neq i_{0}}\overline{\Delta}_{i_{0}j}{\bf x}_{j}+\frac{1}{2}\sum\limits_{i\in I^{-i_{0}}}\sum\limits_{j\neq i, j\in I^{-i_{0}}}{\bf x}_{i}^{\top}\overline{\Delta}_{ij}{\bf x}_{j}\\
&: =f^{(i_{0})}({\bf x})+f^{(-i_{0})}({\bf x}).
\end{aligned}
\end{equation}
We have the following result.
\begin{proposition}\label{prop2}
$f^{(i)}({\bf x})={\bf x}_{i}^{\top}\nabla_{{\bf x}_{i}}f({\bf x})$, $i\in\{1,2,\cdots,n\}$.
\end{proposition}
{\bf Proof.}
According to \eqref{reform-1} and \eqref{reform-2}, we know that 
$$
f({\bf x})=\frac{1}{2}{\bf x}^{\top}\Delta {\bf x}=\sum\limits_{i}\sum\limits_{j>i}{\bf x}_{i}^{\top}\overline{\Delta}_{ij}{\bf x}_{j}.
$$
Then $\nabla_{{\bf x}_{i}}f({\bf x})=\sum\limits_{j\neq i}^{n}\overline{\Delta}_{ij}{\bf x}_{j}$.
Due to \eqref{QP4}, we have
$$
f^{(i_{0})}({\bf x})={\bf x}_{i_{0}}^{\top}\sum\limits_{j\neq i_{0}}\overline{\Delta}_{i_{0}j}{\bf x}_{j}={\bf x}_{i_{0}}^{\top}\nabla_{{\bf x}_{i_{0}}}f({\bf x}).
$$
Then for every index $i\in\{1,2,\cdots,n\}$, $f^{(i)}({\bf x})={\bf x}_{i}^{\top}\nabla_{{\bf x}_{i}}f({\bf x}).$
\hfill$\Box$

\begin{remark}\label{remark}
The above proposition demonstrates that for each block ${\bf x}_{i}$, $f({\bf x})$ is basically a linear function with respect to ${\bf x}_{i}$. This is a key property which will be further explored in the subsequent of the paper.
\end{remark}

\section{Continuous relaxation of \texorpdfstring{$\eqref{QP6}$}{} and two methods}\label{sec-3}
In this section, we relax the sparse constraint problem \eqref{QP6} to a continuous problem and propose two methods to solve the relaxation problem.

\subsection{Relaxation problem}\label{sec-3-1}
Both \eqref{QP2} and \eqref{QP6} are essentially a discrete optimization problem, which is in general NP hard and therefore is extremely difficult to solve. A popular way to deal with \eqref{QP2} or \eqref{QP6} is to relax the discrete constraint and consider solving the relaxed continuous problem. By removing the last constraint in \eqref{QP6}, we obtain the following relaxation problem
\begin{equation}\label{QP-relax}
\begin{aligned}
\min\limits_{{\bf x}\in\mathbb{R}^{nK}}& \hspace{2mm}f({\bf x})\\
{\rm s.t.}&\hspace{2mm}{\bf x}_{i}^{\top}{\bf 1}_{k}=1,\ i=1,\cdots,n,\\
&\hspace{2mm}{\bf x}\geq0.
\end{aligned}
\end{equation}
\eqref{QP-relax} is a continuous problem with simplex constraints. Due to Prop.~\ref{prop2}, the following result holds by Thm. 1 in \cite{no37}, which address the relation between \eqref{QP6} and the relaxation problem \eqref{QP-relax}.
\begin{theorem}\label{thm-1}
There exists a global minimum ${\bf x}^{*}$ of the problem \eqref{QP-relax} such that $\|{\bf x}^{*}\|_{0}=n$, and the global minimizer ${\bf x}^{*}\in\mathbb{R}^{nK}$ is also a global minimizer of the problem \eqref{QP6}.
\end{theorem}

Based on Thm.~\ref{thm-1}, starting from any global minimizer of problem \eqref{QP-relax}, we can eventually find a point ${\bf x}^{*}$ that is a global minimizer of both \eqref{QP-relax} and \eqref{QP6}. Therefore, we can find a global minimizer of \eqref{QP6} using the approach described in Alg.~\ref{alg:1} below.

\begin{algorithm}
\caption{The procedure of finding a global minimizer of \eqref{QP6}}
\label{alg:1}
\begin{algorithmic}[1]
\STATE \textbf{Input}: a global minimizer of \eqref{QP-relax}: ${\bf y}^{0}=(({\bf y}^{0}_{1})^{\top}, ({\bf y}^{0}_{2})^{\top}, \cdots, ({\bf y}^{0}_{n})^{\top})^{\top}\in\mathbb{R}^{nK}$. Let ${\bf x}=(({\bf x}_{1})^{\top}, ({\bf x}_{2})^{\top}, \cdots, ({\bf x}_{n})^{\top})^{\top}={\bf 0}\in\mathbb{R}^{nK}$;
\FOR {$i=1:n$}
\STATE For $i$-th block ${\bf y}^{0}_{i}$ of $y$, we find an index $p^{i}$ which $({\bf y}^{0}_{i})_{p^{i}}\geq({\bf y}^{0}_{i})_{q}$, $q=1,2,\cdots,K$.
\STATE Let $({\bf x}_{i})_{p^{i}}=1$.
\ENDFOR
\STATE \textbf{Output}: ${\bf x}=(({\bf x}_{1})^{\top}, ({\bf x}_{2})^{\top}, \cdots, ({\bf x}_{n})^{\top})^{\top}\in\mathbb{R}^{nK}$, which is a global minimizer of \eqref{QP6}.
\end{algorithmic}
\end{algorithm}

\subsection{Quadratic penalty method}\label{sec-3-2}
To solve the continuous relaxation problem \eqref{QP-relax}, various optimization methods can be used. It has been verified in \cite{no36,no37} that the quadratic penalty method is highly efficient in solving such kind of problem. Therefore, in this paper, we continue to apply the quadratic penalty method to solve \eqref{QP-relax}. The idea of this method is as follows. Due to the result in Thm.~\ref{thm-1}, there is no need to solve \eqref{QP-relax} to get an accurate global minimizer. All we need is to identify the support set of the global minimizer of \eqref{QP-relax} so that we can use Alg.~\ref{alg:1} to obtain a global minimizer of \eqref{QP6}. Therefore, we can penalize the equality constraint and in each iteration, we solve the quadratic penalty problem
\begin{equation}\label{PL}
\min\limits_{0\leq {\bf x}\leq M} g^{l}({\bf x}):=f({\bf x})+\frac{\theta_{l}}{2}\sum\limits_{i=1}^{n}({\bf x}_{i}^{\top}{\bf 1}_{k}-1)^{2},
\end{equation}
where $ M\geq1$ is a given value and the upper bound ${\bf x}\leq M$ is added to make sure that the penalty problem is well-defined. Details of the quadratic penalty method are given in Alg.~\ref{alg:2}.
\begin{algorithm}
\caption{Quadratic penalty method }
\label{alg:2}
\begin{algorithmic}[1]
\STATE \textbf{Input}: an initial point ${\bf x}^{0}$, a parameter $\theta_{0}>0$, $l=0$;
\FOR{$l$}
\STATE Solve \eqref{PL} to get a ${\bf x}^{l}$.
\IF {the termination rule is not satisfied}
\STATE Choose parameter $\theta_{l+1}\geq\theta_{l}$; $l=l+1$;
\ELSE
\STATE break;
\ENDIF
\ENDFOR
\STATE Apply Alg.~\ref{alg:1} to ${\bf x}^{l}$ to get a global minimizer of \eqref{QP6}.
\end{algorithmic}
\end{algorithm}

\begin{theorem}\label{thm-3}{\rm\cite{nno5,nno6}}
Let $\{{\bf x}^{l}\}$ be the sequence generated by Alg.~\ref{alg:2}, and assume that ${\bf x}^{l}$ is a global minimizer of \eqref{PL}. Let $\lim\limits_{l\rightarrow\infty}\theta_{l}=+\infty$. Then any accumulation point of this generated sequence is a global minimum of \eqref{QP-relax}.
\end{theorem}
Thm.~\ref{thm-3} addresses the convergence of the quadratic penalty method. A detailed proof can be found in \cite{nno5,nno6}.

\begin{assumption}\label{assum1}
Let $\{{\bf x}^{l}\}$ be the sequence generated by Alg.~\ref{alg:2}, with $\lim\limits_{l\rightarrow\infty}\theta_{l}=+\infty$. Here, $l$ is a positive integer. Suppose $\lim\limits_{l\rightarrow\infty}{\bf x}^{l}={\bf z}$, and ${\bf z}$ is a global minimizer of \eqref{QP-relax}.
\end{assumption}

Define the support set at ${\bf x}$ by $\Gamma({\bf x}):=\{j:\ {\bf x}_{j}>0\}$. Based on Thm. 3 and Thm. 4 in \cite{no37}, the following result holds.
\begin{theorem}\label{thm-supportset}
Under Assump.~\ref{assum1}, there exists a global minimizer ${\bf z}^{*}$ of \eqref{QP6} such that for $l$ sufficiently large, it holds that $\Gamma({\bf x}^{l})=\Gamma({\bf z}^{*})$.
\end{theorem}

Thm.~\ref{thm-supportset} shows that under Assump.~\ref{assum1}, the support set of the global minimizer for the original problem \eqref{QP6} can be precisely recovered when the number of iterations is sufficiently large.

Inspired by \cite{no37}, we employ a projected Newton method based on the active set to solve the subproblem \eqref{PL}, as detailed in Alg. 3 of \cite{no37}. Our primary goal is to identify the support set of the global minimizer of \eqref{QP6}, rather than its magnitude. Therefore, the projected Newton method is applied to the nonlinear problem \eqref{PL} with simple box constraints. 

\subsection{Quadratic penalty regularized method}\label{sec-3-3}

To solve \eqref{QP-relax}, we propose the following quadratic penalty regularized problem
\begin{equation}\label{PL_BB}
\min\limits_{0\leq {\bf x}\leq M} h^{l}({\bf x}):=f({\bf x})+\frac{\theta_{l}}{2}\sum\limits_{i=1}^{n}({\bf x}_{i}^{\top}{\bf 1}_{k}-1)^{2}+\lambda_{l}\|{\bf x}\|^{2}_{2},
\end{equation}
where $\theta_{l}>0$ and $\lambda_{l}>0$ are both parameters. The motivation is that we know $f({\bf x})$ is not necessarily a convex function due to the properties of $\Delta$. Therefore, to ensure algorithmic robustness, we introduce an $L_{2}$-norm term for ${\bf x}$.
When $\lambda_{l}$ tends to $0$ as $l$ grows, \eqref{PL_BB} exhibits similarly favorable behavior to \eqref{PL}, since \eqref{PL} has been proven to have excellent numerical performance. Details of the quadratic penalty regularized method are given in Alg.~\ref{alg:3}.


\begin{algorithm}
\caption{Quadratic penalty regularized method}
\label{alg:3}
\begin{algorithmic}[1]
\STATE \textbf{Input}: an initial point ${\bf x}^{0}$, max iterations $T$, the parameter $\theta_{0}>0$, $\lambda_{0}>0$, $l=0$;
\FOR{$l$}
\STATE Solve \eqref{PL_BB} to obtain a ${\bf x}^{l}$.
\IF {the termination rule is not satisfied}
\STATE Choose parameter $\theta_{l+1}\geq\theta_{l}$, $\lambda_{l+1}\leq\lambda_{l}$; $l=l+1$;
\ELSE
\STATE break;
\ENDIF
\ENDFOR
\STATE Apply Alg.~\ref{alg:1} to ${\bf x}^{l}$ to get a global minimizer of \eqref{QP6}.
\end{algorithmic}
\end{algorithm}

\begin{theorem}\label{thm-NSPG}
Let $\{{\bf x}^{l}\}$ be the sequence generated by Alg.~\ref{alg:3}. Suppose that ${\bf x}^{l}$ is a global minimizer of \eqref{PL_BB}. Let $\lim\limits_{l\rightarrow\infty}\theta_{l}=+\infty$ and $\lim\limits_{l\rightarrow\infty}\lambda_{l}=0$. Then any accumulation point of the sequence $\{{\bf x}^{l}\}$ is a global minimum of \eqref{QP-relax}. 
\end{theorem}
{\bf Proof.}
Let $\bar{{\bf x}}$ be a global minimizer of \eqref{PL_BB}, that is, $f(\bar{{\bf x}})\leq f({\bf x})$ for all ${\bf x}_{i}$ with $({\bf x}_{i})^{\top}{\bf 1}_{k}-1=0$, $i=1,\cdots,n$. Since ${\bf x}^{l}$ minimizes $h^{l}(\cdot \ ;\theta_{l},\lambda_{l})$ for each $l$, we have $h^{l}({\bf x}^{l};\theta_{l},\lambda_{l})\leq h^{l}(\bar{{\bf x}};\theta_{l},\lambda_{l})$, which leads to 
\begin{equation}\label{thm-NSPG1}
\begin{aligned}
&\hspace{5.3mm}f({\bf x}^{l})+\frac{\theta_{l}}{2}\sum\limits_{i=1}^{n}(({\bf x}^{l}_{i})^{\top}{\bf 1}_{k}-1)^{2}+\lambda_{l}\|{\bf x}^{l}\|^{2}_{2}\\
&\leq f(\bar{{\bf x}})+\frac{\theta_{l}}{2}\sum\limits_{i=1}^{n}((\bar{{\bf x}}_{i})^{\top}{\bf 1}_{k}-1)^{2}+\lambda_{l}\|\bar{{\bf x}}\|^{2}_{2}\\
&=f(\bar{{\bf x}})+\lambda_{l}\|\bar{{\bf x}}\|^{2}_{2}.
\end{aligned}
\end{equation}
Then we have $$\sum\limits_{i=1}^{n}(({\bf x}^{l}_{i})^{\top}{\bf 1}_{k}-1)^{2}\leq \frac{2}{\theta_{l}}\left(f(\bar{{\bf x}})-f({\bf x}^{l})+\lambda_{l}(\|\bar{{\bf x}}\|^{2}_{2}-\|{\bf x}^{l}\|^{2}_{2})\right).$$
Suppose that $\lim\limits_{l\rightarrow\infty}{\bf x}^{l}={\bf z}$, we obtain that
\begin{equation}\notag
\begin{aligned}
\sum\limits_{i=1}^{n}(({\bf z}_{i})^{\top}{\bf 1}_{k}-1)^{2}&=\lim\limits_{l\rightarrow\infty}\sum\limits_{i=1}^{n}(({\bf x}^{l}_{i})^{\top}{\bf 1}_{k}-1)^{2}\\
&\leq \lim\limits_{l\rightarrow\infty}\frac{2}{\theta_{l}}\left(f(\bar{{\bf x}})-f({\bf x}^{l})+\lambda_{l}(\|\bar{{\bf x}}\|^{2}_{2}-\|{\bf x}^{l}\|^{2}_{2})\right)=0,
\end{aligned}
\end{equation}
where $\lim\limits_{l\rightarrow\infty}\theta_{l}=+\infty$. Hence, we have that $(({\bf z}_{i})^{\top}{\bf 1}_{k}-1)^{2}=0$ for all $i=1,\cdots,n$, so that ${\bf z}$ is feasible. Moreover, by $\lim\limits_{l\rightarrow\infty}\lambda_{l}=0$ and taking the limit as $l\rightarrow\infty$ in \eqref{thm-NSPG1}, we have that
\begin{equation}\notag
\begin{aligned}
\lim\limits_{l\rightarrow\infty}\left(f({\bf x}^{l})+\frac{\theta_{l}}{2}\sum\limits_{i=1}^{n}(({\bf x}^{l}_{i})^{\top}{\bf 1}_{k}-1)^{2}+\lambda_{l}\|{\bf x}^{l}\|^{2}_{2}\right)
&\leq \lim\limits_{l\rightarrow\infty}\left(f(\bar{{\bf x}})+\lambda_{l}\|\bar{{\bf x}}\|^{2}_{2}\right)\\
&\leq f(\bar{{\bf x}}).
\end{aligned}
\end{equation}
Then according to nonnegativity of $\theta_{l}$ and of each $(({\bf x}^{l}_{i})^{\top}{\bf 1}_{k}-1)^{2}$, we obtain that
\begin{equation}\notag
\begin{aligned}
f({\bf z})&\leq f({\bf z})+\lim\limits_{l\rightarrow\infty}\left(\frac{\theta_{l}}{2}\sum\limits_{i=1}^{n}(({\bf x}^{l}_{i})^{\top}{\bf 1}_{k}-1)^{2}+\lambda_{l}\|{\bf x}^{l}\|^{2}_{2}\right)\\
&\leq f(\bar{{\bf x}}).
\end{aligned}
\end{equation}
Since ${\bf z}$ is a feasible point whose objective value is no larger than that of the global solution $\bar{{\bf x}}$, we obtain that ${\bf z}$ is a global minimum of \eqref{QP-relax}.
\hfill$\Box$

Thm.~\ref{thm-NSPG} addresses the convergence of the quadratic penalty regularized method. Here, we employ a spectral projected gradient algorithm proposed by \cite{NSPG} to solve \eqref{PL_BB}. 

\section{Numerical Results}\label{sec-4}
In this section, we evaluate the performance of our algorithm on different data sets.
We compare the solution quality of quadratic penalty algorithm (QP-GC) and quadratic penalty regularized method (QPR-GC) with the leading commercial solver Gurobi \cite{Gurobi} for \eqref{QP1} and the Boltzmann machine (BM) from the latest research \cite{GC-BM}. Additionally, we conduct case studies on two real-world graphs. 
To comprehensively demonstrate the superiority of our proposed methods, we include comparisons not only with BM and Gurobi but also with Louvain algorithm \cite{louvain}, currently the most popular clustering approach.
We implement the Boltzmann machine by ourselves, since the original code is not available. Our algorithm, the Boltzmann machine and Gurobi run on a 16-core/8-thread machine. All experiments are written in MATLAB R2023a\footnote{\url{https://ww2.mathworks.cn}} running in Windows 11 on a 12th Gen Intel(R) Core(TM) i7-12800HX CPU at 2.00 GHz with 128 GB of RAM, and all graphs are generated using Python 3.8\footnote{\url{https://www.python.org/downloads/}}. 

\subsection{Evaluation of clustering quality}\label{ecaluation-cq}
To start with, we introduce some measures to evaluate the clustering quality. We use the comparison of intra-cluster density, inter-cluster density, and overall density as metrics to evaluate clustering quality. It is demonstrated in \cite{cq-1} and \cite{cq-2} that the evaluation of the quality of the cluster using these metrics is far superior to the most popular graph cluster quality function, i.e., modularity \cite{cq-3}. 

The number of edges connecting vertices within cluster $i$ is denoted as $|E_{ii}|$, and the number of edges connecting a vertex in cluster $i$ to a vertex in cluster $j$ is denoted as $|E_{ij}|$.
The overall density is defined by
\begin{equation*}
\kappa=\frac{|E(G)|}{0.5\times n(n-1)}.
\end{equation*}
For a cluster $i$, the intra-cluster density is given by 
\begin{equation*}
\kappa^{(i)}_{intra}=\frac{|E_{ii}|}{0.5\times n^{(i)}(n^{(i)}-1)}.
\end{equation*}
For a cluster $i$ and a cluster $j$, $i\neq j$, the inter-cluster density is given by 
\begin{equation*}
\kappa^{(ij)}_{inter}=\frac{|E_{ij}|}{n^{(i)}\times n^{(j)}}.
\end{equation*}
For a graph clustered into $K$ clusters, the mean intra-cluster density and the mean inter-cluster density are given respectively by 
\begin{equation*}
\bar{\kappa}_{intra}=\frac{1}{K}\sum\limits^{K}_{i=1}\kappa^{(i)}_{intra},\ \bar{\kappa}_{inter}=\frac{1}{0.5\times K(K-1)}\sum\limits^{K}_{i=1}\sum\limits^{K}_{j=i+1}\kappa^{(ij)}_{inter}.
\end{equation*}

Using the quantities defined above, we can measure the quality of clustering. A high quality clustering groups vertices into clusters such that, on average, the links between vertices within these clusters are denser than the links between vertices in different clusters \cite{cq-1,cq-2}. Therefore, a high quality clustering should satisfy the following inequality
\begin{equation}\label{cq-standard}
\bar{\kappa}_{inter}<\kappa<\bar{\kappa}_{intra}.
\end{equation}

For the choice of distance, considering factors such as wide applicability, computational complexity, and interpretability, we take the Jaccard distance as the measurement to report our result. Furthermore, it is both demonstrate in \cite{jaccard-4,jaccard-3} that the Jaccard distance is the most suitable distance metric for graph clustering.

\subsection{Synthetic graphs}\label{sec-6-1}
In this part, we perform experiments on fifteen different synthetic graphs, generated using two different models: the Planted Partition Model (PPM) \cite{PPM} and the Stochastic Block Model (SBM) \cite{SBM}. These graphs have known classification properties, as detailed in Tab.~\ref{table-1}. 

We first test on three PPM graph models with $250$ vertices, each containing five clusters with $50$ vertices per cluster. Secondly, we test on six SBM graph models with $3000$ vertices and six with $6000$ vertices. Each graph with $3000$ vertices contains $30$ clusters, where cluster sizes range from $25$ to $200$ vertices. Similarly, each graph with $6000$ vertices comprises $60$ clusters, with cluster sizes varying between $35$ and $200$ vertices. For each intra-cluster edge probability, two different inter-cluster edge probabilities are used to generate the graphs. These graphs are generated using the following intra-cluster/inter-cluster edge probabilities (${\rm P_{intra}}$/${\rm P_{inter}}$) as shown in Tab.~\ref{table-1}.
All synthetic graphs are generated using the NetworkX Python library \cite{no26}. The clustering process becomes more complicated as intra-cluster edge probability decreases, inter-cluster edge probability increases, and cluster sizes vary. An increase in the number of graph vertices also adds to the computational difficulty. 
\begin{table}[ht]
\centering
\caption{Details of fifteen synthetic graphs.}
\label{table-1}
\begin{tabular}{ccccccc}
\hline\hline
\rowcolor[HTML]{FFFFFF} 
Model type & Graph name & ${\rm P_{intra}}$ & ${\rm P_{inter}}$ & $K$  & $n$    & Cluster sizes \\ \hline\hline
\rowcolor[HTML]{EFEFEF} 
PPM        & G1\_PPM    & 0.9          & 0.1          & 5  & 250  & 50            \\ \hline
\rowcolor[HTML]{FFFFFF} 
PPM        & G2\_PPM    & 0.85         & 0.15         & 5  & 250  & 50            \\ \hline
\rowcolor[HTML]{EFEFEF} 
PPM        & G3\_PPM    & 0.8          & 0.2          & 5  & 250  & 50            \\ \hline
\rowcolor[HTML]{FFFFFF} 
SBM        & G1\_SBM3K    & 0.9          & 0.05         & 30 & 3000 & {[}25,200{]}  \\ \hline
\rowcolor[HTML]{EFEFEF} 
SBM        & G2\_SBM3K    & 0.9          & 0.1          & 30 & 3000 & {[}25,200{]}  \\ \hline
\rowcolor[HTML]{FFFFFF} 
SBM        & G3\_SBM3K    & 0.85         & 0.05         & 30 & 3000 & {[}25,200{]}  \\ \hline
\rowcolor[HTML]{EFEFEF} 
SBM        & G4\_SBM3K    & 0.85         & 0.1          & 30 & 3000 & {[}25,200{]}  \\ \hline
\rowcolor[HTML]{FFFFFF} 
SBM        & G5\_SBM3K    & 0.8          & 0.05         & 30 & 3000 & {[}25,200{]}  \\ \hline
\rowcolor[HTML]{EFEFEF} 
SBM        & G6\_SBM3K    & 0.8          & 0.1          & 30 & 3000 & {[}25,200{]}  \\ \hline
\rowcolor[HTML]{FFFFFF} 
SBM        & G1\_SBM6K    & 0.9          & 0.05         & 60 & 6000 & {[}35,200{]}  \\ \hline
\rowcolor[HTML]{EFEFEF} 
SBM        & G2\_SBM6K    & 0.9          & 0.1          & 60 & 6000 & {[}35,200{]}  \\ \hline
\rowcolor[HTML]{FFFFFF} 
SBM        & G3\_SBM6K    & 0.85         & 0.05         & 60 & 6000 & {[}35,200{]}  \\ \hline
\rowcolor[HTML]{EFEFEF} 
SBM        & G4\_SBM6K    & 0.85         & 0.1          & 60 & 6000 & {[}35,200{]}  \\ \hline
\rowcolor[HTML]{FFFFFF} 
SBM        & G5\_SBM6K    & 0.8          & 0.05         & 60 & 6000 & {[}35,200{]}  \\ \hline
\rowcolor[HTML]{EFEFEF} 
SBM        & G6\_SBM6K    & 0.8          & 0.1          & 60 & 6000 & {[}35,200{]}  \\ \hline
\end{tabular}
\end{table}

It is worth noting that the mean intra-cluster density $\bar{\kappa}_{intra}$ (or the mean inter-cluster density $\bar{\kappa}_{inter}$) is an empirical estimate of the intra-cluster edge probability ${\rm P_{intra}}$ (or the inter-cluster edge probability ${\rm P_{inter}}$) in the generated models, as has been demonstrated by \cite{GC-BM}. Therefore, the closer $\bar{\kappa}_{intra}$ (or $\bar{\kappa}_{inter}$) calculated from the clustering obtained by a given method is to ${\rm P_{intra}}$ (or ${\rm P_{inter}}$) of the generation model, the higher the quality of the solution generated by that method. This also indicates that the clustering performance of the method is better. In the results presentation, we use $\epsilon_{\text{intra}}$ (or $\epsilon_{\text{inter}}$) to evaluate clustering quality. Here, $\epsilon_{intra}$ (or $\epsilon_{inter}$) represents the absolute value of the difference between the algorithm-computed $\bar{\kappa}_{intra}$ (or $\bar{\kappa}_{inter}$) and ${\rm P_{intra}}$ (or ${\rm P_{inter}}$) on the graph. We say that the smaller its value, the more accurate the clustering result, especially with regard to $\epsilon_{intra}$.

\subsubsection{Performance for QP-GC and QPR-GC}\label{sec-6-1-1}
\paragraph{Performance for QP-GC}\label{sec-6-1-1-1}
The proposed QP-GC method for solving graph clustering problems is tailored to the problem \eqref{QP1}. This is due to the objective function and its gradient given by $\nabla_{{\bf x}_{i}}f({\bf x})=\sum_{j\neq i}\overline{\Delta}_{ij}{\bf x}_{j}$. Additionally, \eqref{QP1} involves only a single vector variable ${\bf x}>0$ and $n$ linear equality constraints, all of which are linear.  For Hessian computation, we use sparsity to reduce computational complexity from $O(n^{2}k^{2})$ to $O(nk^{2})$. These characteristics ensure that the computational process runs efficiently.

The parameters for QP-GC are set as follows. In Alg.~\ref{alg:2}, we set the initial point ${\bf x}^{0}$ as a random matrix. Update $\theta_{l}$ as
\begin{equation*}
\theta_{l+1}=\left\{\begin{array}{ll}
		3\theta_{l},&{\rm if} \ \|h^{l}\|>0.01 \ {\rm and} \ \theta_{l}<\bar{\theta},\\
		\theta_{l},&{\rm otherwise},
		\end{array}	\right.
\end{equation*}
where $h^{l}=\left(({\bf x}^{l}_{1})^{\top}{\bf 1}-1, \cdots, ({\bf x}^{l}_{n})^{\top}{\bf 1}-1\right)$. For the experiments on the PPM model, we set $\bar{\theta}=10^{5}$ and $\theta_{0}=10$, while for the experiments on the SBM model, we set $\bar{\theta}=10^{10}$ and $\theta_{0}=2500$. 
Each ${\bf x}^{l}$ returned is projected onto a binary assignment matrix using Alg.~\ref{alg:1}. For the experiments on the small-scale PPM model, we set $\epsilon=0.0009$, while for the large-scale SBM model experiments, we set $\epsilon=2\times 10^{-8}$. For more complex large-scale cases, extensive experiments demonstrate that setting termination criteria by limiting the support set size enables early iteration termination, effectively reducing unnecessary computational overhead. Specifically, if the size of the support set found in Alg. 3 of \cite{no37} is greater than or equal to $90\%$ to $95\%$ of n×k in the graph, the iteration is terminated.

\begin{table}[ht]
\centering
\caption{Numerical results of QP-GC for three PPM models with $n=250$ and $K=5$.}
\setlength{\tabcolsep}{4.5pt} 
\resizebox{\linewidth}{!}{
\begin{tabular}{cccccccccc}
    \toprule
    \multirow{2.4}{*}{Graph name} & \multicolumn{3}{c}{Graph characteristics} & \multicolumn{6}{c}{QP-GC}\\
    \cmidrule(r){2-4} \cmidrule(r){5-10}   
    & ${\rm P_{intra}}$ & ${\rm P_{inter}}$ & $\kappa$ & $\bar{\kappa}_{intra}$ & $\bar{\kappa}_{inter}$ & $\epsilon_{intra}$ & $\epsilon_{inter}$ & times(s) & $C_{iter}$ \\
    \midrule
    G1\_PPM & 0.9   & 0.1   & 0.26 & 0.90 & 0.10 & 0 & 0 & 0.4 & 220 \\
    \midrule
    G2\_PPM & 0.85  & 0.15  & 0.29 & 0.85 & 0.15 & 0 & 0 & 0.5 & 340 \\
    \midrule
    G3\_PPM & 0.8   & 0.2   & 0.32 & 0.80 & 0.20 & 0 & 0 & 0.7 & 430 \\
    \bottomrule
    \end{tabular}%
    }
  \label{tab:2}%
\end{table}%

\begin{table}[htbp]
\centering
\caption{Numerical results of QP-GC for six SBM models with $n=3000$ and $K=30$. All times are formatted as mm:ss (minutes:seconds).}
\setlength{\tabcolsep}{4.5pt} 
\resizebox{\linewidth}{!}{
\begin{tabular}{cccccccccc}
    \toprule
    \multirow{2.4}{*}{Graph name} & \multicolumn{3}{c}{Graph characteristics} & \multicolumn{6}{c}{QP-GC}\\
    \cmidrule(r){2-4} \cmidrule(r){5-10}   
    & ${\rm P_{intra}}$ & ${\rm P_{inter}}$ & $\kappa$ & $\bar{\kappa}_{intra}$ & $\bar{\kappa}_{inter}$ & $\epsilon_{intra}$ & $\epsilon_{inter}$ & times & $C_{iter}$ \\
    \midrule
    G1\_SBM3K & 0.90  & 0.05  & 0.09 & 0.78 & 0.06 & 0.12 & 0.01 & 13:29 & 2960 \\
    \midrule
    G2\_SBM3K & 0.90  & 0.10  & 0.14 & 0.73 & 0.12 & 0.17 & 0.02 & 30:52 & 6000 \\
    \midrule
    G3\_SBM3K & 0.85  & 0.05  & 0.09 & 0.73 & 0.06 & 0.12 & 0.01 & 21:24 & 4530 \\ 
    \midrule
    G4\_SBM3K & 0.85  & 0.10 & 0.13 & 0.71 & 0.11 &  0.14 & 0.01  & 16:39 & 3790 \\
    \midrule
    G5\_SBM3K & 0.80 & 0.05 & 0.08 & 0.65 & 0.06 & 0.15 & 0.01 & 22:42 & 4970 \\
    \midrule
    G6\_SBM3K & 0.80 & 0.10 & 0.13 & 0.66 & 0.11 & 0.14 & 0.01 & 33:52 & 6660 \\
    \bottomrule
    \end{tabular}%
    }
  \label{tab:3}%
\end{table}%

\begin{table}[htbp]
\centering
\caption{Numerical results of QP-GC for six SBM models with $n=6000$ and $K=60$. All times are formatted as hh:mm:ss (hours:minutes:seconds).}
\setlength{\tabcolsep}{4.2pt} 
\resizebox{\linewidth}{!}{
\begin{tabular}{cccccccccc}
    \toprule
    \multirow{2.4}{*}{Graph name} & \multicolumn{3}{c}{Graph characteristics} & \multicolumn{6}{c}{QP-GC}\\
    \cmidrule(r){2-4} \cmidrule(r){5-10}   
    & ${\rm P_{intra}}$ & ${\rm P_{inter}}$ & $\kappa$ & $\bar{\kappa}_{intra}$ & $\bar{\kappa}_{inter}$ & $\epsilon_{intra}$ & $\epsilon_{inter}$ & times & $C_{iter}$ \\
    \midrule
    G1\_SBM6K & 0.90  & 0.05  & 0.09 & 0.74 & 0.06 & 0.16 & 0.01 & 1:54:49 & 5920 \\
    \midrule
    G2\_SBM6K & 0.90  & 0.10  & 0.14 & 0.73 & 0.11 & 0.17 & 0.01 & 2:19:00 & 7160 \\
    \midrule
    G3\_SBM6K & 0.85  & 0.05  & 0.09 & 0.69 & 0.06 & 0.21 & 0.01 & 2:24:59 & 6000 \\ 
    \midrule
    G4\_SBM6K & 0.85  & 0.10 & 0.13 & 0.68 & 0.11 &  0.22 & 0.01  & 2:27:36 & 7400 \\
    \midrule
    G5\_SBM6K & 0.80 & 0.05 & 0.08 & 0.64 & 0.06 & 0.26 & 0.01 & 2:02:24 & 6180 \\
    \midrule
    G6\_SBM6K & 0.80 & 0.10 & 0.13 & 0.64 & 0.11 & 0.26 & 0.01 & 2:31:20 & 7440 \\
    \bottomrule
    \end{tabular}%
    }
  \label{tab:3-6000}%
\end{table}%

Tab.~\ref{tab:2}, Tab.~\ref{tab:3} and Tab.~\ref{tab:3-6000} show the clustering results of QP-GC for three PPM models, six SBM models with $3000$ vertices, and six more complex SBM models with $6000$ vertices, respectively. We also report the CPU time and the number of iterations $C_{iter}$ for QP-GC. Our findings show that for PPM models, the algorithm achieves perfect clustering results in under one second. For the six SBM models with $3000$ vertices, it obtains reasonably good results within one hour, while the more complex cases with $6000$ vertices require approximately two hours.

\paragraph{Performance for QPR-GC}\label{sec-6-1-1-2}
QPR-GC method is also specifically designed for graph clustering problems. As shown in Alg.~\ref{alg:3}, we employ the BB1 step size \cite{BB3} here, with an initial step size of $0.1$. For PPM models, we set the penalty parameter $\theta_{l}$ as a fixed value of $300$. The parameter $\lambda_{l}$ starts at $1$ and decays uniformly in logarithmic space to $10^{-7}$ over $100$ iterations, then remains fixed at $10^{-7}$ until termination. Additionally, we set $\eta=0.5$, $\gamma=0.5$, and $tol=10^{-6}$. 
For more complex SBM models, the penalty parameter $\theta_{l}$ takes a fixed value within $[1000,3000]$. $\lambda_{l}$ is initialized at $0.1$ and decays logarithmically to $10^{-8}$ within $100$ iterations, maintaining this value thereafter. Other parameters are set as $\eta=0.6$, $\gamma=0.8$, and $tol=10^{-8}$. For this problem, we additionally establish the following termination criterion: the iterative process terminates and returns ${\bf x}^{l}$ once any vector ${\bf x}^{l+1}_{i}$ $(i = 1,\cdots,n)$ is identified as the zero vector.

\begin{table}[ht]
\centering
\caption{Numerical results of QPR-GC for three PPM models with $n=250$ and $K=5$.}
\setlength{\tabcolsep}{4.4pt} 
\resizebox{\linewidth}{!}{
\begin{tabular}{cccccccccc}
    \toprule
    \multirow{2.4}{*}{Graph name} & \multicolumn{3}{c}{Graph characteristics} & \multicolumn{6}{c}{QPR-GC}\\
    \cmidrule(r){2-4} \cmidrule(r){5-10}   
    & ${\rm P_{intra}}$ & ${\rm P_{inter}}$ & $\kappa$ & $\bar{\kappa}_{intra}$ & $\bar{\kappa}_{inter}$ & $\epsilon_{intra}$ & $\epsilon_{inter}$ & times(s) & $C_{iter}$ \\
    \midrule
    G1\_PPM & 0.9   & 0.1   & 0.26 & 0.90 & 0.10 & 0 & 0 & 0.1 & 108 \\
    \midrule
    G2\_PPM & 0.85  & 0.15  & 0.29 & 0.85 & 0.15 & 0 & 0 & 0.2 & 158 \\
    \midrule
    G3\_PPM & 0.8   & 0.2   & 0.32 & 0.80 & 0.20 & 0 & 0 & 0.2 & 165 \\
    \bottomrule
    \end{tabular}%
    }
  \label{tab:2-BB}%
\end{table}%

\begin{table}[htbp]
\centering
\caption{Numerical results of QPR-GC for six SBM models with $n=3000$ and $K=30$. All times are formatted as mm:ss (minutes:seconds).}
\setlength{\tabcolsep}{4.5pt} 
\resizebox{\linewidth}{!}{
\begin{tabular}{cccccccccc}
    \toprule
    \multirow{2.4}{*}{Graph name} & \multicolumn{3}{c}{Graph characteristics} & \multicolumn{6}{c}{QPR-GC}\\
    \cmidrule(r){2-4} \cmidrule(r){5-10}   
    & ${\rm P_{intra}}$ & ${\rm P_{inter}}$ & $\kappa$ & $\bar{\kappa}_{intra}$ & $\bar{\kappa}_{inter}$ & $\epsilon_{intra}$ & $\epsilon_{inter}$ & times& $C_{iter}$ \\
    \midrule
    G1\_SBM3K & 0.90  & 0.05  & 0.09 & 0.79 & 0.07 & 0.11 & 0.02 & 10:39 & 3500 \\
    \midrule
    G2\_SBM3K & 0.90  & 0.10  & 0.14 & 0.73 & 0.12 & 0.17 & 0.02 & 24:54 & 7000 \\
    \midrule
    G3\_SBM3K & 0.85  & 0.05  & 0.09 & 0.73 & 0.06 & 0.12 & 0.01 & 18:46 & 6000 \\ 
    \midrule
    G4\_SBM3K & 0.85  & 0.10 & 0.13 & 0.71 & 0.11 & 0.14 & 0.01  & 15:57 & 5000 \\
    \midrule
    G5\_SBM3K & 0.80 & 0.05 & 0.08 & 0.65 & 0.06 & 0.15 & 0.01 & 20:52 & 7000 \\
    \midrule
    G6\_SBM3K & 0.80 & 0.10 & 0.13 & 0.65 & 0.11 & 0.15 & 0.01 & 22:05 & 6800 \\
    \bottomrule
    \end{tabular}%
    }
  \label{tab:3-BB}%
\end{table}%

\begin{table}[htbp]
\centering
\caption{Numerical results of QPR-GC for six SBM models with $n=6000$ and $K=60$. All times are formatted as hh:mm:ss (hours:minutes:seconds).}
\setlength{\tabcolsep}{4.2pt} 
\resizebox{\linewidth}{!}{
\begin{tabular}{cccccccccc}
    \toprule
    \multirow{2.4}{*}{Graph name} & \multicolumn{3}{c}{Graph characteristics} & \multicolumn{6}{c}{QPR-GC}\\
    \cmidrule(r){2-4} \cmidrule(r){5-10}   
    & ${\rm P_{intra}}$ & ${\rm P_{inter}}$ & $\kappa$ & $\bar{\kappa}_{intra}$ & $\bar{\kappa}_{inter}$ & $\epsilon_{intra}$ & $\epsilon_{inter}$ & times & $C_{iter}$ \\
    \midrule
    G1\_SBM6K & 0.90  & 0.05  & 0.09 & 0.72 & 0.06 & 0.18 & 0.01 & 2:09:35 & 10000 \\
    \midrule
    G2\_SBM6K & 0.90  & 0.10  & 0.14 & 0.65 & 0.11 & 0.25 & 0.01 & 2:13:51 & 10000 \\
    \midrule
    G3\_SBM6K & 0.85  & 0.05  & 0.09 & 0.65 & 0.06 & 0.20 & 0.01 & 2:19:16 & 10000 \\ 
    \midrule
    G4\_SBM6K & 0.85  & 0.10 & 0.13 & 0.56 & 0.11 & 0.29 & 0.01  & 2:00:37 & 9000 \\
    \midrule
    G5\_SBM6K & 0.80 & 0.05 & 0.08 & 0.55 & 0.06 & 0.25 & 0.01 & 2:23:23 & 10000 \\
    \midrule
    G6\_SBM6K & 0.80 & 0.10 & 0.13 & 0.49 & 0.11 & 0.31 & 0.01 & 2:11:27 & 10000 \\
    \bottomrule
    \end{tabular}%
    }
  \label{tab:3-BB6000}%
\end{table}%

Following the same presentation approach as for the QP-GC algorithm results, we present Tab.~\ref{tab:2-BB}, Tab.~\ref{tab:3-BB} and Tab.~\ref{tab:3-BB6000} to demonstrate the clustering results of the QPR-GC algorithm for three PPM models, six SBM models with $3000$ vertices, and six more complex SBM models with $6000$ vertices, respectively. We also report the CPU time of QPR-GC and the iteration counts $C_{iter}$. The results show that for PPM models, perfectly accurate clustering solutions are obtained within one second. For the six SBM models with $3000$ vertices, reasonably good results are achieved within half an hour, while approximately two hours are required for the more complex cases with $6000$ vertices.

\subsubsection{Comparison between QP-GC and QPR-GC}\label{sec-6-1-2}
\paragraph{Performance comparison of algorithms}\label{sec-6-1-2-1}
We know that QP-GC is a second-order method based on the active-set projected Newton approach, which requires the computation of the Hessian matrix. In contrast, QPR-GC is a first-order algorithm based on the BB step sizes and nonmonotone line search, and only requires the computation of the gradient. Each iteration of QP-GC involves greater computational load, primarily from PCG (preconditioned conjugate gradient) solving and Hessian construction. While QPR-GC has lower per iteration costs, it may require more iterations to converge. A detailed comparison of computational complexity per iteration for QP-GC and QPR-GC is provided in Tab.~\ref{tab:compare1}.

\begin{table}[htbp]
  \centering
  \renewcommand{\arraystretch}{1.3}
  \caption{Comparison of computational complexity between QP-GC and QPR-GC.}
  \setlength{\tabcolsep}{4.3pt} 
  \resizebox{\linewidth}{!}{
    \begin{tabular}{c|c|c}
    \Xhline{1.2pt}
    Per-iteration steps   & QP-GC  & QPR-GC  \\
    \Xhline{1.2pt}
    Gradient computation & $O(kn^{2})$ & $O(kn^{2})$ \\
    \hline
    Hessian computation & $O(nk^{2})$ & Not required \\
    \hline
    Linear system solving &  PCG solving $O(nk^{2})$ & Not required \\
    \hline
    \multirow{2}{*}{Step size selection } & \multirow{2}{*}{Armijo line search $O(kn^{2})$} & \multirow{2}{*}{\shortstack{BB step size $O(kn)$ and \\ Non-monotonic line search $O(kn^{2})$}} \\
          &       &  \\
    \hline
    Others    & Active set identification $O(kn)$ & Projection operation $O(kn)$ \\
    \hline
    Total    & $O(kn^{2}+nk^{2})$ & $O(kn^{2})$ \\
    \Xhline{1.2pt}
    \end{tabular}%
    }
  \label{tab:compare1}%
\end{table}%

\begin{figure}[ht]
	\centering
	\begin{minipage}{1\linewidth}	
		\subfigure[G1\_PPM]{
			\label{fig:a1QPPG}
			\includegraphics[width=0.32\linewidth,height=1.4in]{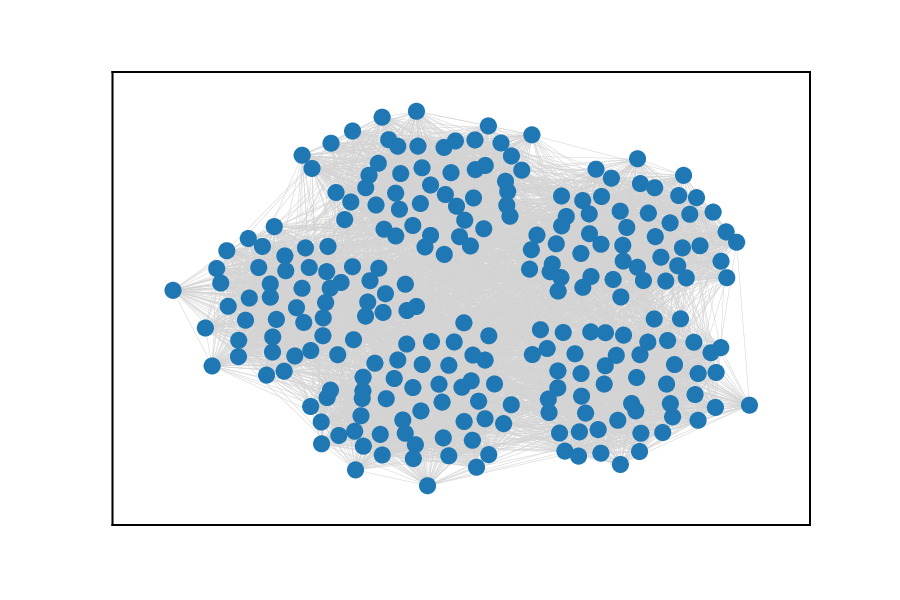}	
		}\hspace{-0.45em}
		\subfigure[QP-GC clustering result]{
			\label{fig:b1_ppm_QPPG}
			\includegraphics[width=0.32\linewidth,height=1.4in]{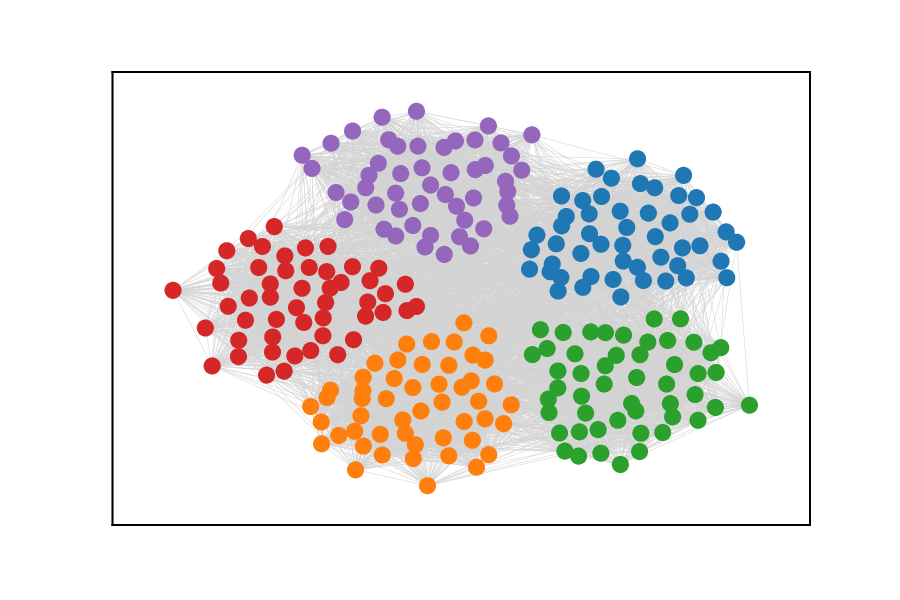}
		}\hspace{-0.45em}
		\subfigure[QPR-GC clustering result]{
			\label{fig:b1_ppm_BB}
			\includegraphics[width=0.32\linewidth,height=1.4in]{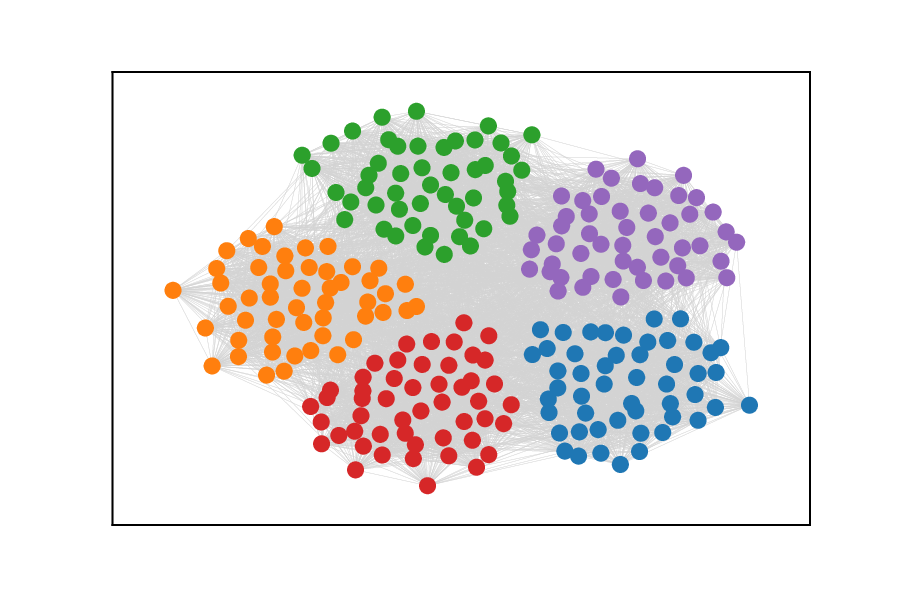}
		}
	\end{minipage}
    	 \vskip -0.3cm 
	\begin{minipage}{1\linewidth}
		\subfigure[G2\_PPM]{
			\label{fig:a2QPPG}
			\includegraphics[width=0.32\linewidth,height=1.4in]{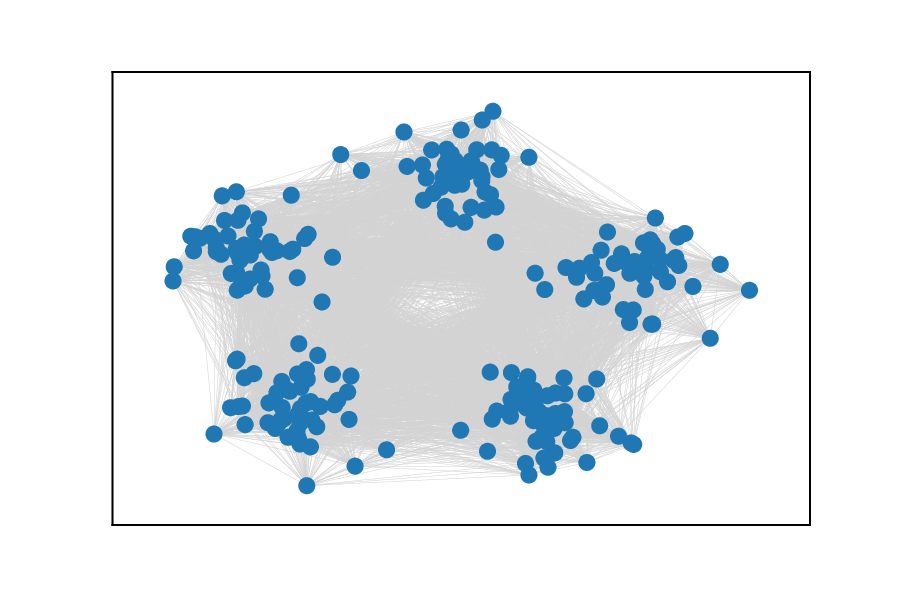}		
		}\hspace{-0.45em}
		\subfigure[QP-GC clustering result]{
			\label{fig:b2_ppm_QPPG}
			\includegraphics[width=0.32\linewidth,height=1.4in]{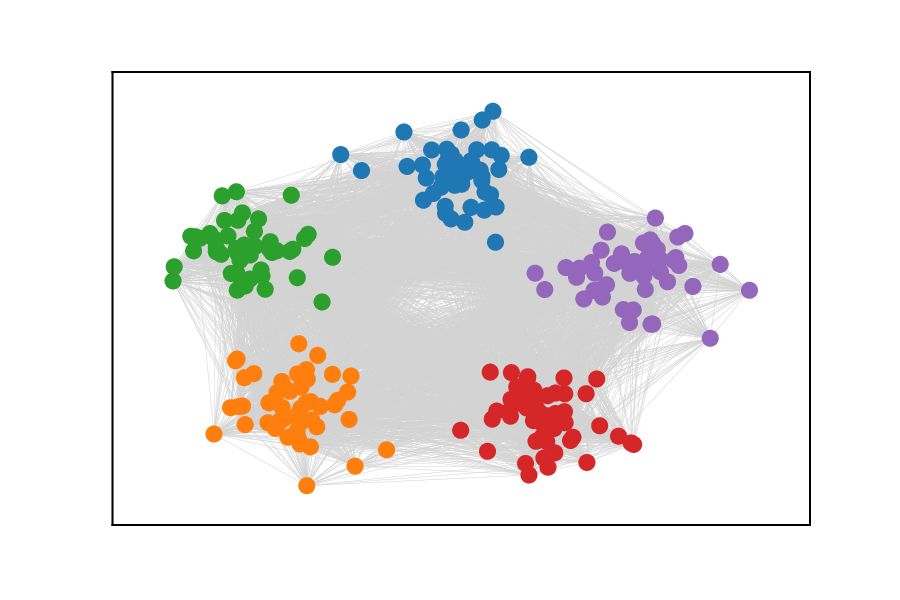}
		}\hspace{-0.45em}
		\subfigure[QPR-GC clustering result]{
			\label{fig:b2_ppm_BB}
			\includegraphics[width=0.32\linewidth,height=1.4in]{figure/ppm2_QPPG.eps}
		}
	\end{minipage}
       	\vskip -0.3cm 
	\begin{minipage}{1\linewidth}
		\subfigure[G3\_PPM]{
			\label{fig:a3QPPG}
			\includegraphics[width=0.32\linewidth,height=1.4in]{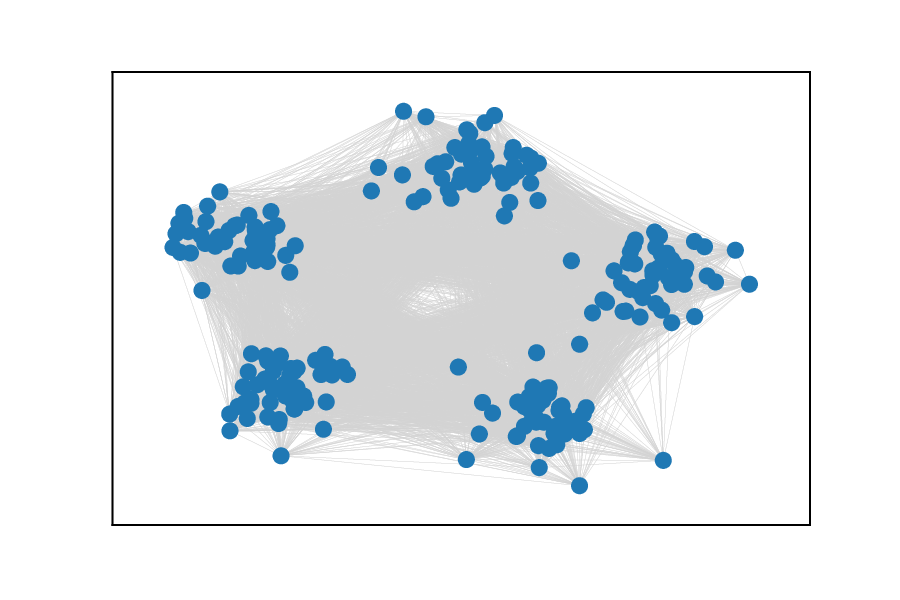}		
		}\hspace{-0.45em}
		\subfigure[QP-GC clustering result]{
			\label{fig:b3_ppm_QPPG}
			\includegraphics[width=0.32\linewidth,height=1.4in]{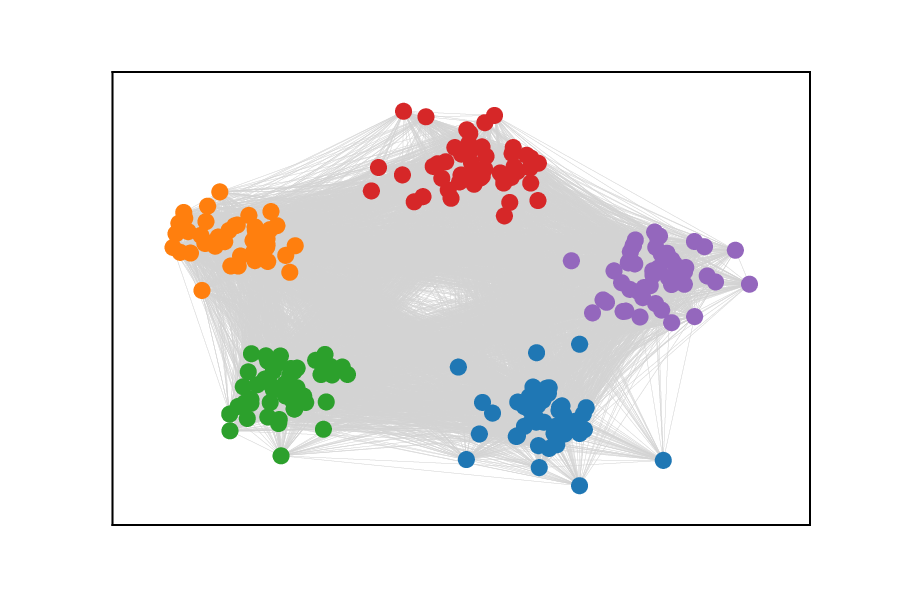}
		}\hspace{-0.45em}
		\subfigure[QPR-GC clustering result]{
			\label{fig:b3_ppm_BB}
			\includegraphics[width=0.32\linewidth,height=1.4in]{figure/ppm3_QPPG.eps}
		}
	\end{minipage}
	\caption{A visual representation of three PPM graph models after clustering with our method }
	\label{fig-2}
\end{figure}

\paragraph{Comparison of Clustering Results}\label{sec-6-1-2-2}
To provide a more intuitive demonstration of QP-GC and QPR-GC, we visualize the above three PPM graphs before and after clustering, as shown in Fig.~\ref{fig-2}. The visualizations of the three PPM graphs are shown in the left column of Fig.~\ref{fig-2}. 
The middle column presents the clustering results obtained by QP-GC, where each graph is partitioned into five clusters (represented by distinct colors). Similarly, the right column shows the QPR-GC clustering results, also demonstrating five color partitions. Fig.~\ref{fig-2} illustrates that both methods accurately divide all three graphs into five clusters, indicating high-quality clustering. This visual representation effectively demonstrates the efficacy of QP-GC and QPR-GC.

\begin{figure}[ht]
	\centering
	\begin{minipage}{1\linewidth}	
		\subfigure[the SBM with $3000$ vertices]{
			\label{fig-compare-sbm3000-BBQP}
			\includegraphics[width=0.98\linewidth,height=2.3in]{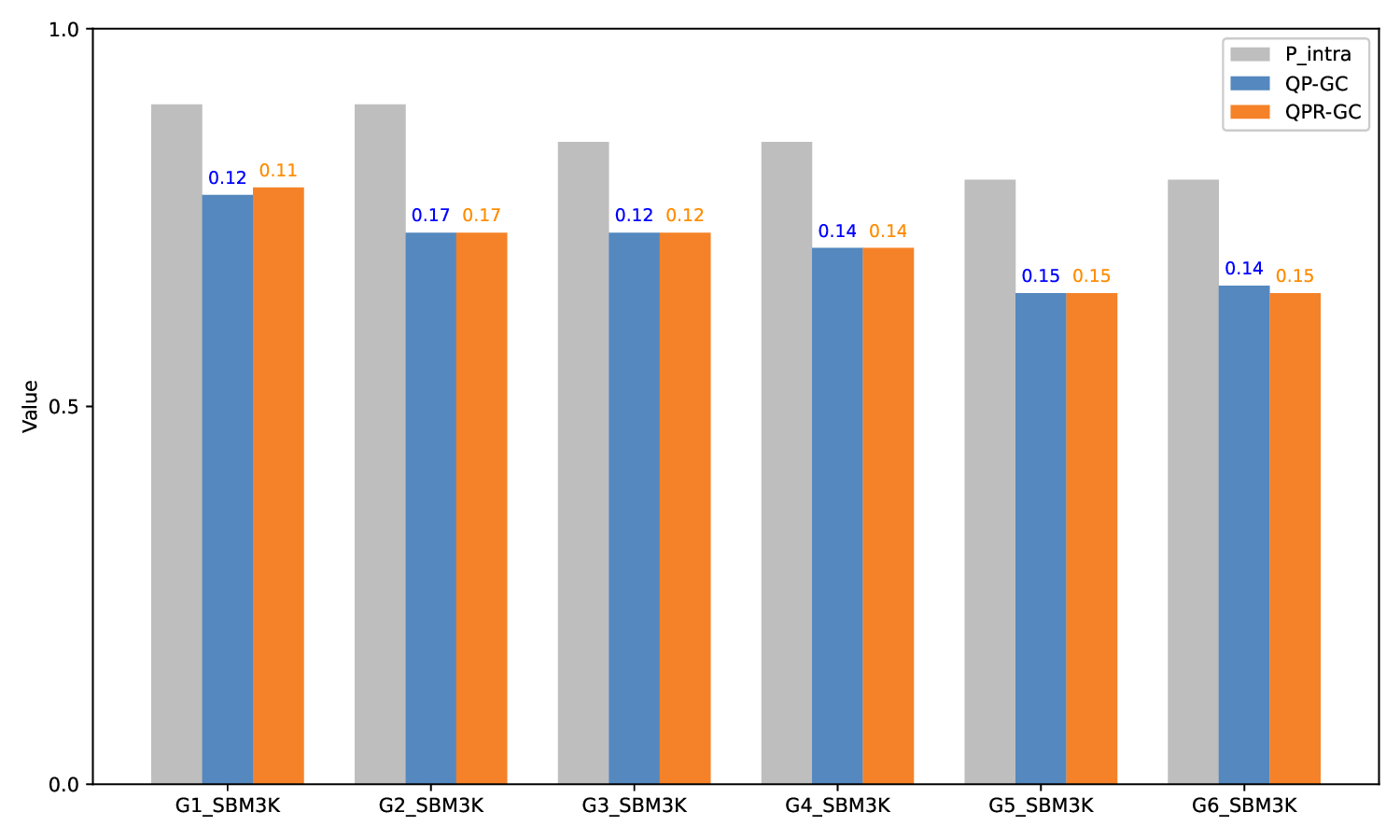}
		}
	\end{minipage}
	\begin{minipage}{1\linewidth}
		\subfigure[the SBM with $6000$ vertices]{
			\label{fig-compare-sbm6000-BBQP}
			\includegraphics[width=0.98\linewidth,height=2.3in]{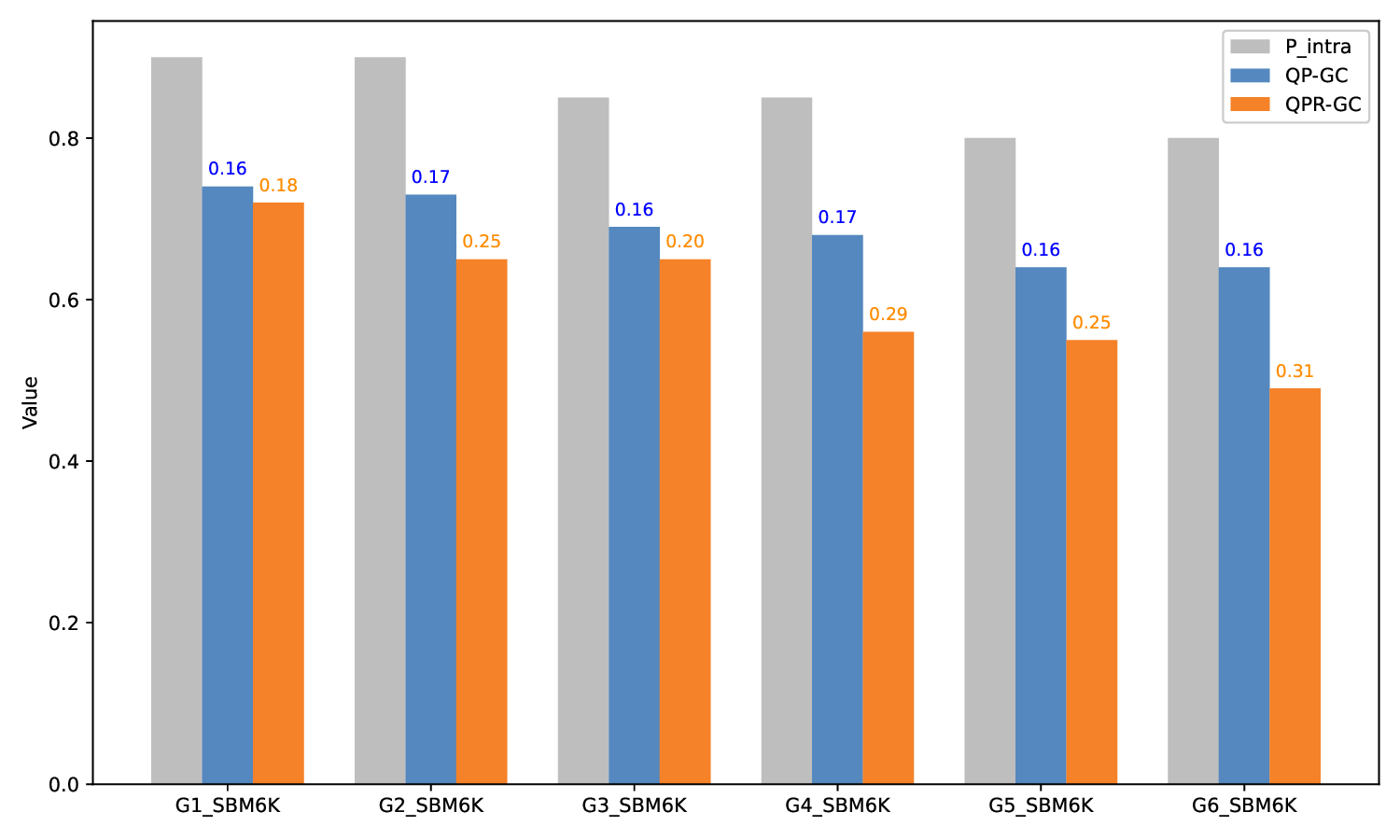}		
		}
	\end{minipage}
	\caption{Dendrogram visualization of intra cluster density distributions: ground truth ${\rm P_{intra}}$ versus $\bar{\kappa}_{intra}$ by QP-GC and QPR-GC for twelve SBM models}
	\label{fig-compare-sbm-BBQP}
\end{figure}

For the more complex set of twelve SBM models, the visualization approach shown in Fig.~\ref{fig-2} proves less intuitive. However, the metric $\epsilon_{intra}$ (or $\epsilon_{inter}$) sufficiently demonstrates clustering performance. We further visualize these values using dendrograms, as presented in Fig.~\ref{fig-compare-sbm-BBQP}. 
The gray bars represent the graph's ${\rm P_{intra}}$, the blue bars show $\bar{\kappa}_{intra}$ obtained by QP-GC, and the orange bars indicate $\bar{\kappa}_{intra}$ computed by QPR-GC. The numerical values displayed above the blue and orange bars correspond to $\epsilon_{intra}$.
From Fig.~\ref{fig-compare-sbm-BBQP}, we observe that for the SBM with $3000$ vertices, there is little difference between QP-GC and QPR-GC. However, for the SBM with $6000$ vertices, QP-GC produces slightly better clustering results than QPR-GC.

\paragraph{Summary}\label{sec-6-1-2-3} 
Based on the above results, we observe that for the PPM and the SBM model with $3000$ vertices, QP-GC and QPR-GC achieve very similar clustering accuracy, but QPR-GC uses less time and requires fewer iterations. However, for the SBM model with $6000$ vertices, QPR-GC performs worse than QP-GC in both clustering accuracy and the CPU time. These findings demonstrate that QP-GC is more suitable for complex, large-scale problems of this type.

\subsubsection{Comparison with other algorithms}\label{sec-6-1-3}

\paragraph{Results on PPM graph models}\label{sec-6-1-3-1}
We compare QP-GC and QPR-GC with Gurobi and BM in terms of $\epsilon_{intra}$, $\epsilon_{inter}$, and CPU time.
From Tab.~\ref{tab:compare-PPM-2}, we can observe that all three solutions produce the same results, and all satisfy \eqref{cq-standard}. In all three cases (G1\_PPM,  G2\_PPM and G3\_PPM), they perfectly recover the generation model. In this example, it can be noticed that under the same conditions, QPR-GC achieves accurate clustering results in a shorter CPU time.

\begin{table}[ht]
  \centering
  \caption{Comparison of some clustering metrics on PPM models with $n=250$ and $K=5$.}
  \setlength{\tabcolsep}{4.5pt} 
  \resizebox{\linewidth}{!}{
  \begin{tabular}{cccccccccc}
    \toprule
    \multirow{2.4}{*}{Methods} & \multicolumn{3}{c}{G1\_PPM} & \multicolumn{3}{c}{G2\_PPM} & \multicolumn{3}{c}{G3\_PPM} \\
    \cmidrule(r){2-4} \cmidrule(r){5-7} \cmidrule(r){8-10} 
    & $\epsilon_{intra}$ & $\epsilon_{inter}$ & times(s)     & $\epsilon_{intra}$ & $\epsilon_{inter}$ & times(s) & $\epsilon_{intra}$ & $\epsilon_{inter}$ & times(s) \\
    \midrule
    QP-GC & \bf{0} & \bf{0}  & 0.4  & \bf{0} & \bf{0} & 0.5 & \bf{0} & \bf{0} & 0.7 \\
    \midrule
    QPR-GC & \bf{0} & \bf{0}  & \bf{0.1}  & \bf{0} & \bf{0} & \bf{0.2} & \bf{0} & \bf{0} & \bf{0.2} \\
    \midrule
    Gurobi & \bf{0} & \bf{0} & 3.2 & \bf{0} & \bf{0} & 3.1 & \bf{0} & \bf{0} & 3.1 \\
    \midrule
    BM & \bf{0}  & \bf{0} & 3.5 & \bf{0} & \bf{0} & 3.6  & \bf{0} & \bf{0} & 3.7 \\
    \bottomrule
    \end{tabular}%
    }
  \label{tab:compare-PPM-2}%
\end{table}%

\begin{figure}[ht]
	\centering
		\subfigure[the $3000$-vertex SBM]{
			\label{fig:intra}
			\includegraphics[width=0.48\linewidth]{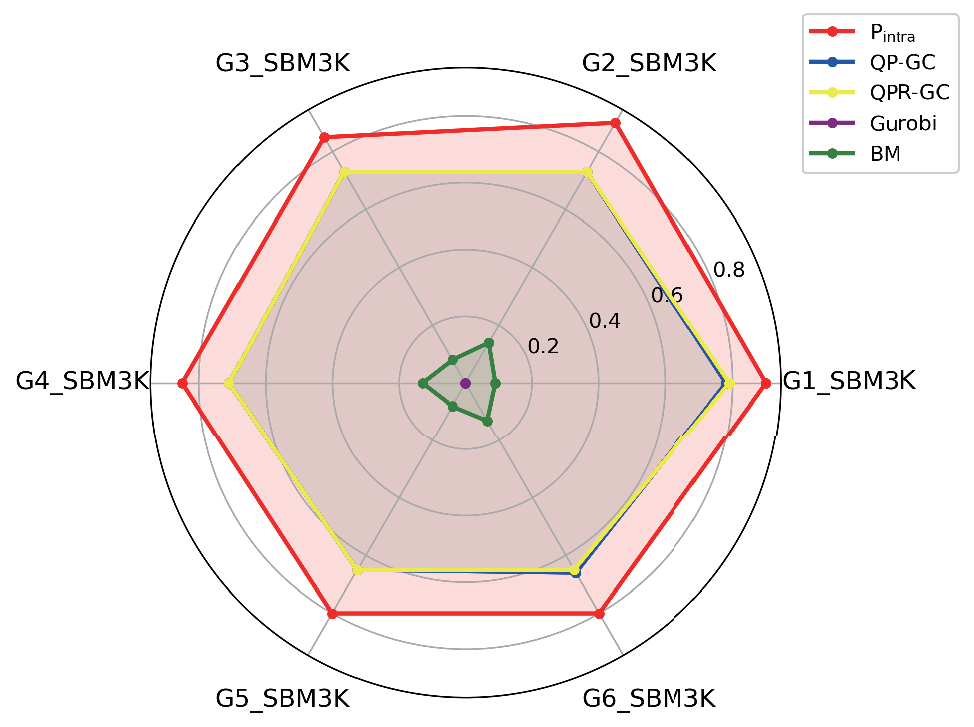}	
		}\noindent
		\subfigure[the $6000$-vertex SBM]{
			\label{fig:intra6000}
			\includegraphics[width=0.48\linewidth]{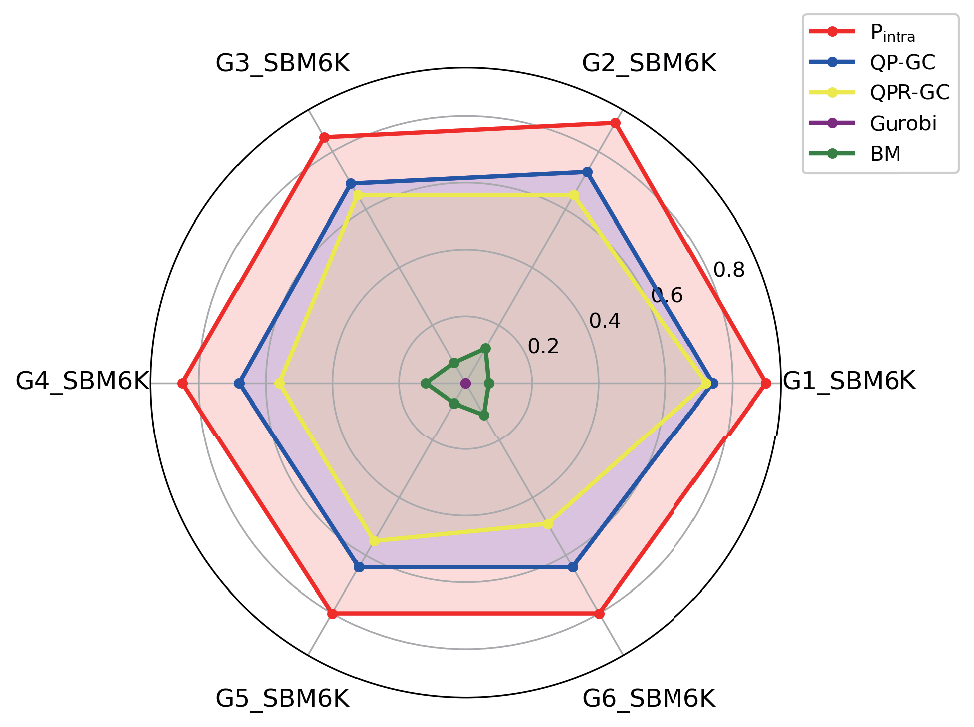}
		}
	\vspace{-0.11in}	
	\caption{The radar plot representation of ${\rm P_{intra}}$ and $\bar{\kappa}_{intra}$ by QP-GC, QPR-GC, Gurobi and BM for twelve SBM models}
	\label{fig:radar}
\end{figure}

\paragraph{Results on SBM graph models}\label{sec-6-1-3-2}
The graph characteristics and numerical results of Gurobi and BM are reported in Tab.~\ref{tab:compare-sbm-1}. We further compare QP-GC and QPR-GC against Gurobi and BM in terms of $\epsilon_{intra}$. To visually demonstrate the superiority of our approach, we present the ${P_{intra}}$ and $\bar{\kappa}_{intra}$ metrics from Tab.~\ref{tab:compare-sbm-1} using radar charts: Fig.~\ref{fig:intra} displays results for the $3000$-vertex SBM, while Fig.~\ref{fig:intra6000} shows the $6000$-vertex case.
Based on the results in Tab.~\ref{tab:compare-sbm-1}, Fig.~\ref{fig:intra}, and Fig.~\ref{fig:intra6000}, we summarize two key observations: 
\begin{itemize}
    \item[(a)] From Tab.~\ref{tab:compare-sbm-1}, we note that BM’s results satisfy inequality \eqref{cq-standard}. However, Gurobi fails to return meaningful results, and the inequality does not hold.
    \item[(b)] From Figs.~\ref{fig:intra} and Figs.~\ref{fig:intra6000}, we observe that regardless of the SBM graph model, QP-GC (blue line) and QPR-GC (yellow line) yield $\bar{\kappa}_{intra}$ values closer to the graph model’s $P_{intra}$ (red line) compared to Gurobi (purple line) and BM (green line).
\end{itemize}
Thus, for larger and more complex SBM graphs, QP-GC and QPR-GC provide faster and more reliable results for these models.

\begin{table}[htbp]
  \centering
  \caption{Numerical results of Gurobi and BM on SBM models: six graphs with $3000$ vertices and six graphs with $6000$ vertices. All times are formatted as hh:mm:ss (hours:minutes:seconds).}
  \setlength{\tabcolsep}{4.3pt} 
  \resizebox{\linewidth}{!}{
    \begin{tabular}{cccccccccc}
    \toprule
    \multirow{2.4}{*}{Graph name} & \multicolumn{3}{c}{Graph characteristics}& \multicolumn{3}{c}{Gurobi} & \multicolumn{3}{c}{BM} \\
    \cmidrule(r){2-4} \cmidrule(r){5-7} \cmidrule(r){8-10}   
    & ${\rm P_{intra}}$ & ${\rm P_{inter}}$ & $\kappa$ & $\bar{\kappa}_{intra}$ & $\bar{\kappa}_{inter}$ & times & $\bar{\kappa}_{intra}$ & $\bar{\kappa}_{inter}$ & times \\
    \midrule
    G1\_SBM3K & 0.90   & 0.05  & 0.09 & 0.00 & 0.01 & 0:13:29 & 0.09 & 0.10 & 0:13:30\\
    \midrule
    G2\_SBM3K & 0.90 & 0.10 & 0.14 & 0.00 & 0.01 & 0:30:52 & 0.14  & 0.15 & 0:30:52 \\
    \midrule
    G3\_SBM3K & 0.85  & 0.05  & 0.09 & 0.00 & 0.01 & 0:21:24 & 0.08 & 0.10 & 0:21:24 \\ 
    \midrule
    G4\_SBM3K & 0.85  & 0.10 & 0.13 & 0.00 & 0.01 & 0:16:39 & 0.13 &  0.14 & 0:16:39 \\
    \midrule
    G5\_SBM3K & 0.80   & 0.05  &  0.08 & 0.00 & 0.01 & 0:22:42 & 0.08 &  0.10 &  0:22:42 \\
    \midrule
    G6\_SBM3K & 0.80 & 0.10  & 0.13 & 0.00 & 0.01 & 0:33:52 & 0.13 &  0.14 & 0:33:52 \\
    \midrule
    G1\_SBM6K & 0.90   & 0.05  & 0.07 & 0.00 & 0.00 & 1:54:49 & 0.07 & 0.07 & 1:54:51\\
    \midrule
    G2\_SBM6K & 0.90 & 0.10 & 0.12 & 0.00 & 0.01 & 2:19:00 & 0.12  & 0.12 & 2:19:02\\
    \midrule
    G3\_SBM6K & 0.85  & 0.05  & 0.07 & 0.00 & 0.00 & 2:24:59 & 0.07 & 0.07 & 2:25:02 \\ 
    \midrule
    G4\_SBM6K & 0.85  & 0.10 & 0.12 & 0.00 & 0.01 & 2:27:36 & 0.12 &  0.12 & 2:27:37 \\
    \midrule
    G5\_SBM6K & 0.80   & 0.05  &  0.07 & 0.00 & 0.00 & 2:02:24 & 0.07 & 0.07 &  2:02:25 \\
    \midrule
    G6\_SBM6K & 0.80 & 0.10  & 0.11 & 0.00 & 0.01 & 2:31:20 & 0.11 &  0.11 & 2:31:22\\
    \bottomrule
    \end{tabular}%
    }
  \label{tab:compare-sbm-1}%
\end{table}%

\subsection{Real-world graphs}\label{sec-6-2}
Although synthetic networks provide a repeatable and controlled testing platform for our experiments, testing algorithms on real-world network data is also necessary, even if real-world graphs are often not the best benchmarks for assessing clustering quality. In this section, we select two representative real-world network datasets: {\it Zachary's Karate Club} (Karate\_club) \cite{no3} and the {\it United States College Football Division IA 2000 season} graph (US\_football\_2000) \cite{no3}. Both contain a known clustering structure. The graph data used in our experiments are available at \url{https://github.com/qinyuenlp/CommunityDetection/tree/master/data}. The characteristics of these graphs are shown in Tab.~\ref{real-graph-c}. Here, we also compare these four methods with the widely popular Louvain algorithm \cite{louvain}, whose code is available at \url{https://github.com/GenLouvain/GenLouvain}.

\begin{table}[ht]
  \centering
  \caption{Graph characteristics of two real-world graphs.}
    \begin{tabular}{ccccc}
    \toprule
    \multirow{2}[4]{*}{Graph name} & \multicolumn{4}{c}{Graph characteristics} \\
\cmidrule{2-5}          & $|V|$   & $|E|$   & $\kappa$ & Clusters  \\
    \midrule
    Karate\_club & 34    & 78    & 0.14  & 2 \\
    \midrule
    US\_football\_2000 & 115   & 613   & 0.09  & 12 \\
    \bottomrule
    \end{tabular}%
  \label{real-graph-c}%
\end{table}%

For experiments on real-world graphs, the following two points must be noted: (1) Real-world graphs are often instances of unknown underlying generative models with random noise. In many cases, modifying the clustering assignment of real-world graph vertices may result in higher intra-cluster density. This means that the true clustering results of the real-world graph vertices may not correspond to the clusters with the highest density. (2) When comparing the clustering results with the true clusters in real-world graphs, it is not sufficient to simply compare the cluster labels. Instead, the comparison should focus on the composition of the clusters. For instance, the cluster $C_{1}$ in the true clustering might not correspond to $C_{1}$ from the QP-GC clustering result, but it might fully correspond to $C_{2}$ from the QP-GC result. These two factors determine that we need to set evaluation criteria different from those used in the synthetic graph experiments.

For both of the following experiments, the QP-GC parameters are set as follows: $\theta_{0} = 20$ and the initial point $S^{0}$ is set as a random matrix. We set $\epsilon=10^{-5}$ in Alg. 3 of \cite{no37}. For the termination condition, we set $\tau=0.01$. Other parameters are the same as the PPM experiments. The QPR-GC parameters are set as follows: (1) For Karate\_club, set $\theta_l = 250$, $m = 5$, $\eta=0.5$, $\gamma=0.5$, and $tol=10^{-7}$; (2) For US\_football\_2000, set $\theta_l = 220$, $m = 28$, $\eta=0.6$, $\gamma=0.8$, and $tol=10^{-6}$.

\subsubsection{Karate club study of Zachary}
\begin{figure}[ht]
	\centering
	\begin{minipage}{1\linewidth}	
		\subfigure[The true clustering]{
			\label{fig:atrue}
			\includegraphics[width=0.47\linewidth,height=1.5in]{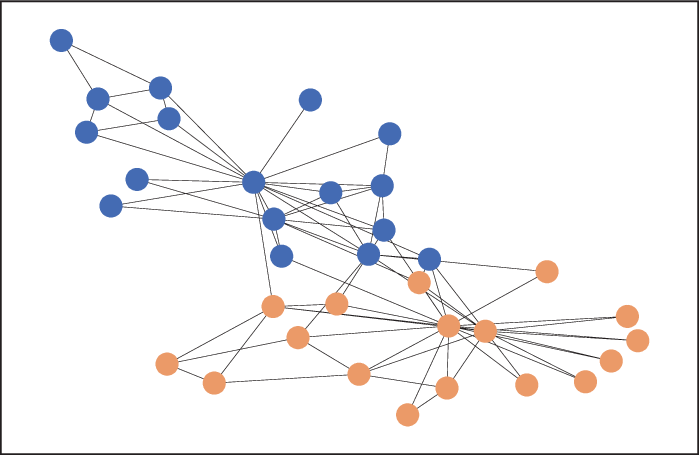}	
		}\noindent
		\subfigure[The clustering result from QP-GC ]{
			\label{fig:bQPPG}
			\includegraphics[width=0.47\linewidth,height=1.5in]{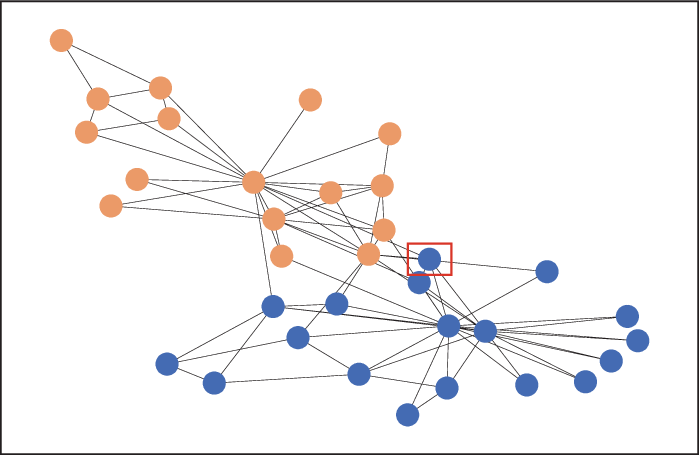}
		}
	\end{minipage}
    	\vskip -0.2cm 
	\begin{minipage}{1\linewidth }
		\subfigure[The clustering result from QPR-GC]{
			\label{fig:cBB}
			\includegraphics[width=0.47\linewidth,height=1.5in]{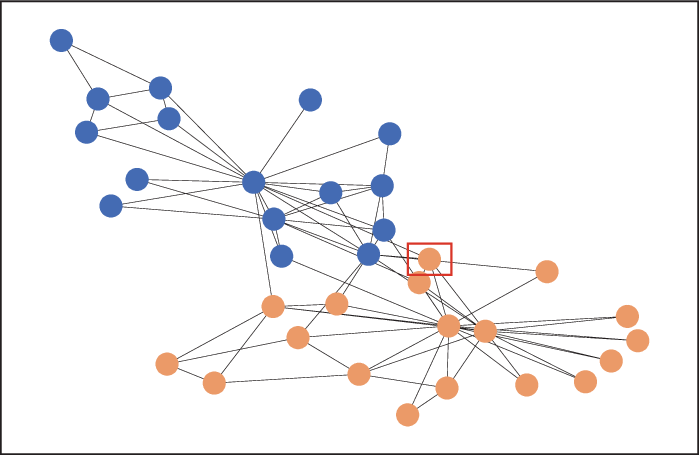}
		}\noindent
		\subfigure[The clustering result from louvain]{
			\label{fig:dlouvain}
			\includegraphics[width=0.47\linewidth,height=1.5in]{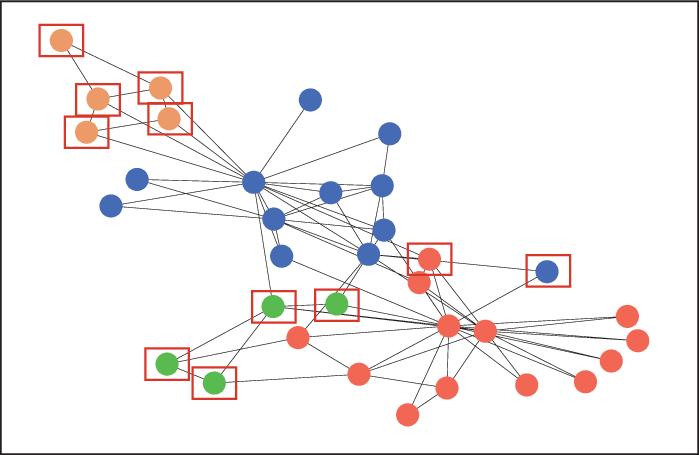}
		}
	\end{minipage}
	\vskip -0.2cm 
	\begin{minipage}{1\linewidth }
		\subfigure[The clustering result from Gurobi]{
			\label{fig:egurobi}
			\includegraphics[width=0.47\linewidth,height=1.5in]{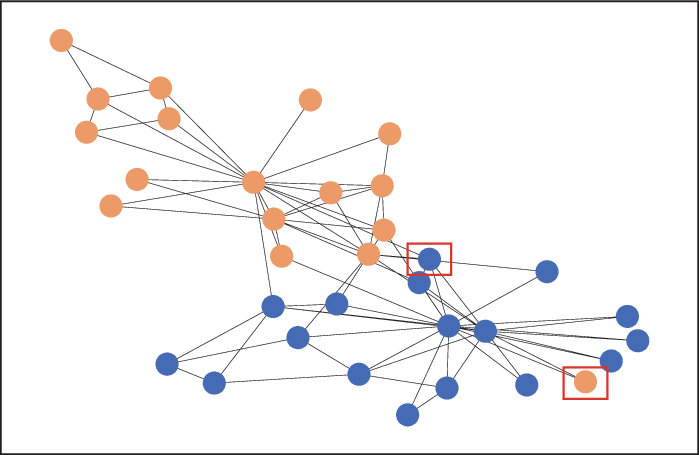}
		}\noindent
		\subfigure[The clustering result from BM]{
			\label{fig:fBM}
			\includegraphics[width=0.47\linewidth,height=1.5in]{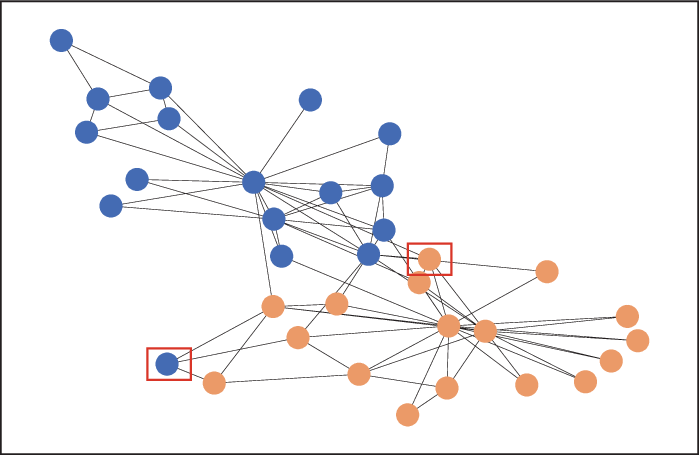}
		}
	\end{minipage}
	\vspace{-0.11in}	
	\caption{The friendship network derived from Zachary's karate club, as described in this paper. Nodes of different colors represent membership in different factions. The red squares in the figure indicate the nodes of clustering error}
	\label{fig:1234}
\end{figure}

The real-world graph network data in this subsection comes from Zachary's famous study of the karate club \cite{no3}. The vertices of the graph represent the $34$ members of the karate club, and the edges represent the social connections between the members. Here, we use the simplified unweighted version of the network. The club split into two factions due to internal conflict, and this division is reflected in the graph structure, as shown in Fig.~\ref{fig:atrue}. Note that in each subfigure of Fig.~\ref{fig:1234}, nodes with the same color indicate that they belong to the same cluster. However, the same color across different subfigures does not convey any particular meaning.

We apply our method to this graph to identify the factions involved in the club's split. In addition, under the same CPU time (approximately $1$ second), we compare the clustering results obtained using Louvain, Gurobi, and BM with those generated by the QP-GC and QPR-GC.
Both Gurobi and BM have been shown to be effective for clustering in \cite{GC-BM}. Since the graph has two cluster structures, we visualize the clustering, making the results visually clear, as shown in Fig.~\ref{fig:1234}. We find that both QP-GC and QPR-GC misclassify only one node, as shown in Fig.~\ref{fig:bQPPG} and Fig.~\ref{fig:cBB}, while Gurobi and BM misclassify two nodes, as shown in Fig.~\ref{fig:egurobi} and Fig.~\ref{fig:fBM}. In contrast, Louvain produces an incorrect four-cluster division, specifically misclassifying $11$ nodes as shown in Figure~\ref{fig:dlouvain}.
Thus, QP-GC and QPR-GC more accurately restore the clustering. 

\subsubsection{US College Football Division IA 2000 season graph}
This subsection uses the well-known United States College Football Division IA $2000$ season graph \cite{no3}, which represents the matchups between $115$ college football teams during the regular season. The vertices in the graph represent the teams, labeled by their university names, and the edges indicate that the two connected teams have faced each other at least once during the regular season. These teams are divided into $12$ conferences, with each conference containing approximately $8$-$12$ teams. Teams within the same conference compete more frequently against each other than against teams from different conferences, resulting in more shared connections among teams in the same conference compared to those in different ones. 

\begin{figure}[ht]
	\centering
	\begin{minipage}{1\linewidth}	
		\subfigure[The true clustering]{
			\label{fig:atrue1}
			\includegraphics[width=0.47\linewidth,height=1.5in]{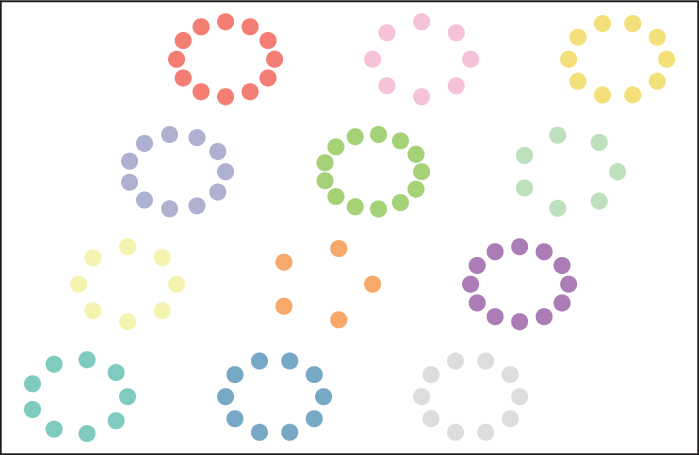}	
		}\noindent
		\subfigure[The clustering result from QP-GC ]{
			\label{fig:bQPPG1}
			\includegraphics[width=0.47\linewidth,height=1.5in]{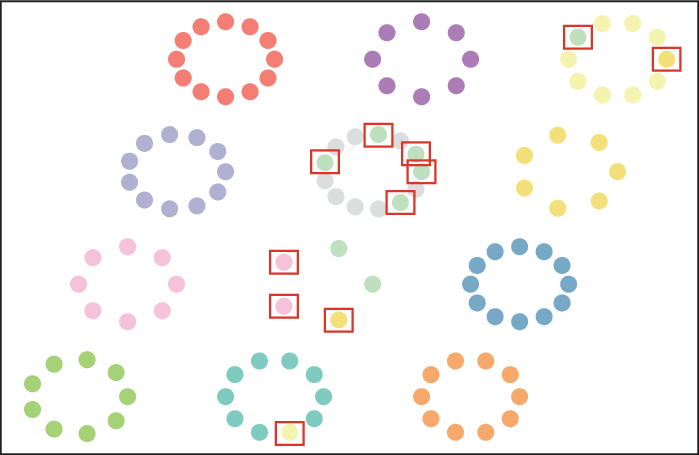}
		}
	\end{minipage}
    	\vskip -0.2cm 
	\begin{minipage}{1\linewidth }
		\subfigure[The clustering result from QPR-GC]{
			\label{fig:cBB1}
			\includegraphics[width=0.47\linewidth,height=1.5in]{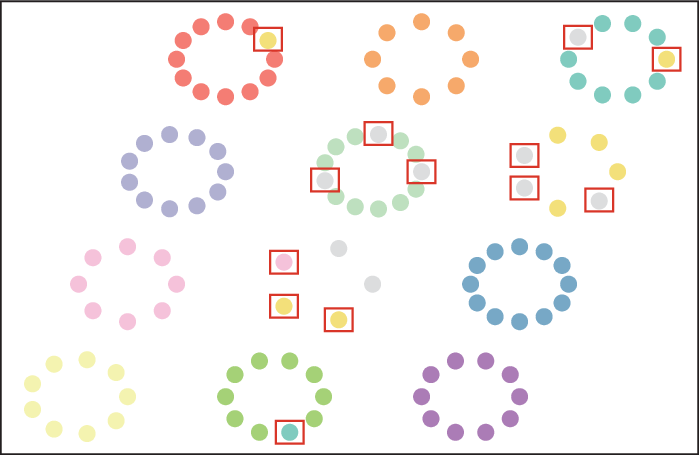}
		}\noindent
		\subfigure[The clustering result from Louvain]{
			\label{fig:dlouvain1}
			\includegraphics[width=0.47\linewidth,height=1.5in]{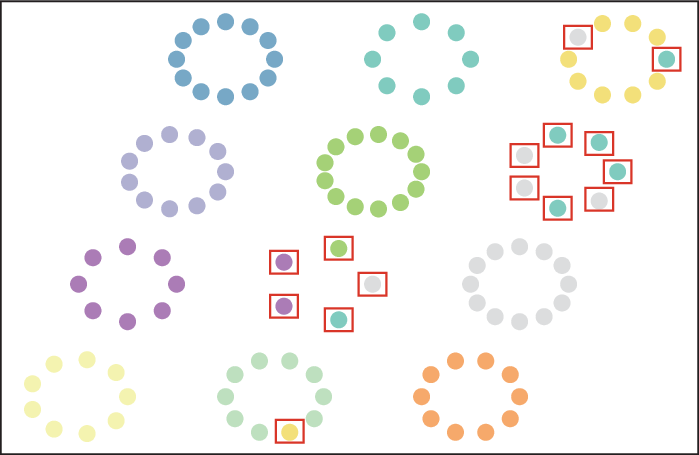}
		}
	\end{minipage}
	\vskip -0.2cm 
	\begin{minipage}{1\linewidth }
		\subfigure[The clustering result from Gurobi]{
			\label{fig:egurobi1}
			\includegraphics[width=0.47\linewidth,height=1.5in]{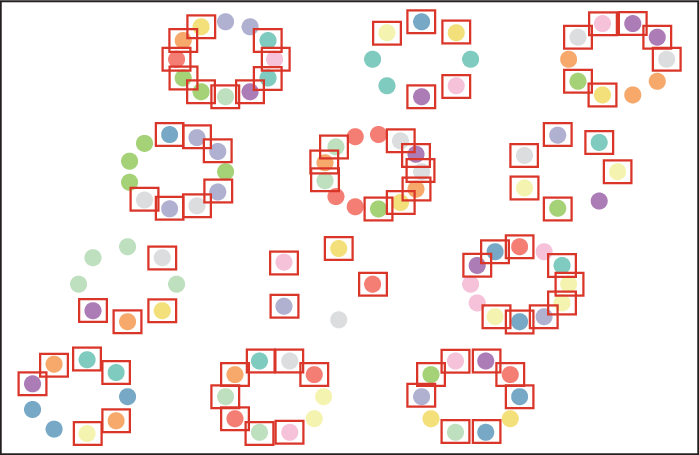}
		}\noindent
		\subfigure[The clustering result from BM]{
			\label{fig:fBM1}
			\includegraphics[width=0.47\linewidth,height=1.5in]{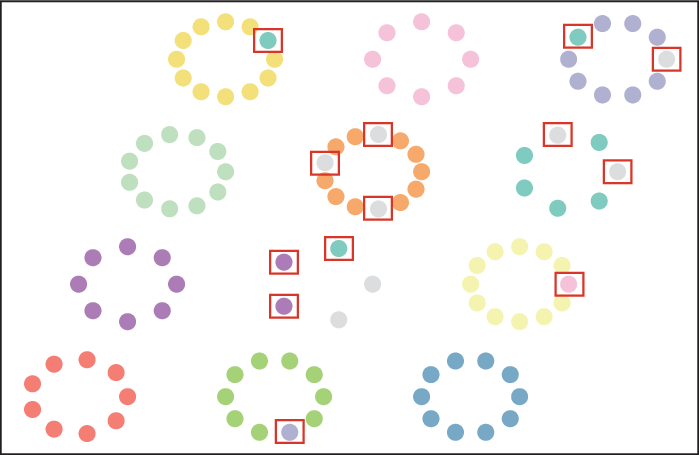}
		}
	\end{minipage}
	\caption{A visualization of the US College Football Division IA 2000 season graph. Each circle represents nodes (teams) belonging to the same cluster, while red squares indicate misclassified nodes}
	\label{fig:football1}
\end{figure}

We also compare the clustering results of QP-GC and QPR-GC with Louvain, Gurobi, and BM when running for the same CPU time (approximately $1$ second). The graph consists of $12$ clusters, as shown in Fig.~\ref{fig:football1}, where each circle represents a cluster. We assign $12$ different colors to these circles. QP-GC misclassifies $11$ nodes, as shown in Fig.~\ref{fig:bQPPG1}, while both QPR-GC and BM misclassify $13$ nodes, as shown in Fig.~\ref{fig:cBB1} and Fig.~\ref{fig:fBM1}. Gurobi performs significantly worse with $83$ misclassified nodes (Fig.~\ref{fig:egurobi1}). Louvain identifies only $10$ clusters (versus the ground truth of $12$) and misclassifies $15$ nodes, as shown in Fig.~\ref{fig:dlouvain1}.
A numerical summary is provided in Tab.~\ref{compare-football-tab}.
Therefore, QP-GC more accurately recovers the clustering structure. 

\begin{table}[ht]
  \centering
  \caption{Numerical Comparison of Clustering Results for US\_football\_2000.}
  \setlength{\tabcolsep}{4.5pt} 
    \begin{tabular}{ccccccc}
    \toprule
    Reference values & \textcolor{blue}{True} & QP-PG & QPR-GC & Louvain & Gurobi & BM \\
    \midrule
   Number of clusters  & \textcolor{blue}{12}  & {\bf 12} &{\bf 12} & 10  & {\bf 12}  & {\bf 12} \\
    \midrule
   Number of misclassified nodes  & \textcolor{blue}{0}  & {\bf 11}   & 13 & 15 & 83 & 13\\
    \bottomrule
    \end{tabular}%
  \label{compare-football-tab}%
\end{table}%

To conclude, we find that QP-GC and QPR-GC demonstrate the ability to recover clusters not only in synthetic graphs but also in real-world graphs.

\section{Conclusion}\label{sec-5}
In this paper, we considered a class of optimization models for graph clustering problems in a unified manner. We reformulated these models as sparse-constrained optimization problems and recovered the solutions of the original problem from the solutions of the relaxed problem without sparse constraints. We apply both the quadratic penalty method and the quadratic penalty regularized method to the relaxation problem. The corresponding subproblems are solved using an active-set projected Newton method and a spectral projected gradient method, respectively. We compare QP-GC and QPR-GC against Louvain, Gurobi, and Boltzmann machines.
Through extensive numerical experiments, we demonstrate that both QP-GC and QPR-GC successfully cluster both synthetic and real-world graphs. For small-scale graphs, both methods achieve satisfactory clustering results, with QPR-GC requiring less computation time. For more complex large-scale graphs, QP-GC not only outperforms other methods in speed but also delivers more accurate clustering results.

\vspace{3.5mm}
\noindent\textbf{Acknowledgements}\hspace{1mm} Not applicable.
\vspace{3.5mm}

\noindent\textbf{Authors contributions}\hspace{1mm} L.Q. reported the optimization model. T.W. proposed the algorithm and conducted numerical experiments. T.W. and L.Q. wrote the main manuscript text. All authors reviewed the manuscript.
\vspace{3.5mm}

\noindent\textbf{Funding}\hspace{1mm} This work is supported by the NSFC 12071032 and NSFC 12271526.
\vspace{3.5mm}

\noindent\textbf{Availability of data and materials}\hspace{1mm} Synthetic graph datasets are provided when needed. The real-world network datasets are public, and the URL of the link to the datasets is provided within the manuscript.

\section*{Declarations}

\noindent\textbf{Consent for publication}\hspace{1mm} All authors consent to the publication of this manuscript.
\vspace{3.5mm}

\noindent\textbf{Competing interests}\hspace{1mm} I declare that the authors have no competing interests as defined by Springer, or other interests that might be perceived to influence the results and/or discussion reported in this paper.
\vspace{3.5mm}

\noindent\textbf{Ethics approval and consent to participate}\hspace{1mm} This research did not involve human participants or animals, and therefore ethical approval and consent were not applicable.



\begin{thebibliography}{9}
\setlength{\itemindent}{-2em}
\setlength{\leftskip}{2em}

\bibitem[Rossi et al. 2020]{no1} Rossi R A, Jin D, Kim S et al (2020) On proximity and structural role-based embeddings in networks: Misconceptions, techniques, and applications. Acm T Knowl Discov D 14(5):1-37. \url{https://doi.org/10.1145/3397191}


\bibitem[Newman 2001]{no2} Newman M E J (2001) The structure of scientific collaboration networks. P Natl Acad Sci 98(2):404-409. \url{https://doi.org/10.1073/pnas.98.2.404}


\bibitem[Girvan et al. 2002]{no3} Girvan M, Newman M E J (2002) Community structure in social and biological networks. P Natl Acad Sci 99(12):7821-7826. \url{https://doi.org/10.1073/pnas.122653799}

\bibitem[Williams et al. 2000]{no4} Williams R J, Martinez N D (2000) Simple rules yield complex food webs. Nature 404(6774):180-183. \url{https://doi.org/10.1038/35004572}

\bibitem[Krause et al. 2003]{no5} Krause A E, Frank K A, Mason D M et al (2003) Compartments revealed in food-web structure. Nature 426(6964):282-285. \url{https://doi.org/10.1038/nature02115}

\bibitem[Jeong et al. 2000]{no6} Jeong H, Tombor B, Albert R (2000) The large-scale organization of metabolic networks. Nature 407(6804):651-654. \url{https://doi.org/10.1038/35036627}

\bibitem[Faloutsos et al. 1999]{no7} Faloutsos M, Faloutsos P, Faloutsos C (1999) On power-law relationships of the internet topology. Comput Commun Rev 29(4):251-262. \url{https://doi.org/10.1145/316194.316229}

\bibitem[Newman 2004a]{no8} Newman M E J (2004) Coauthorship networks and patterns of scientific collaboration. P Natl Acad Sci 101(suppl\_1):5200-5205. \url{https://doi.org/10.1073/pnas.0307545100}

\bibitem[Papadopoulos et al. 2012]{no9} Papadopoulos S, Kompatsiaris Y, Vakali A et al (2012) Community detection in social media: Performance and application considerations. Data Min Knowl Disc 24:515-554. \url{https://doi.org/10.1007/s10618-011-0224-z}

\bibitem[Rostami et al. 2023]{no10} Rostami M, Oussalah M, Berahmand K (2023) Community detection algorithms in healthcare applications: a systematic review. IEEE Access 11:30247-30272. \url{https://doi.org/10.1109/ACCESS.2023.3260652}

\bibitem[Kumar et al. 1999]{no11} Kumar R, Raghavan P, Rajagopalan S (1999) Trawling the web for emerging cyber-communities. Comput Netw 31(11-16):1481-1493. \url{https://doi.org/10.1016/S1389-1286(99)00040-7}

\bibitem[Flake et al. 2000]{no12} Flake G W, Lawrence S, Giles C L (2000) Efficient identification of web communities. In proceedings of the sixth ACM SIGKDD international conference on Knowledge discovery and data mining 2000:150-160. \url{https://doi.org/10.1145/347090.347121}

\bibitem[Steenstrup 2001]{no13} Steenstrup M (2001) Cluster-based networks. In: Ad hoc networking, Addison-Wesley Longman Publishing Co., Inc., USA, 75-138. ISBN 0201309769

\bibitem[Wu et al. 2004]{no14} Wu A Y, Garland M, Han J (2004) Mining scale-free networks using geodesic clustering. In proceedings of the tenth ACM SIGKDD international conference on Knowledge discovery and data mining 2004:719-724. \url{https://doi.org/10.1145/1014052.1014146}

\bibitem[Schaeffer 2007]{no15} Schaeffer S E (2007) Survey: Graph clustering. Comput Sci Rev 1(1):27-64. \url{https://doi.org/10.1016/j.cosrev.2007.05.001}

\bibitem[Hastie et al. 2009]{no16} Hastie T, Tibshirani R, Friedman J (2009) The Elements of Statistical Learning, Second Edition: Data Mining, Inference, and Prediction. Springer, New York. ISBN 978-0-387-84857-0

\bibitem[Miasnikof et al. 2024]{GC-BM} Miasnikof P, Bagherbeik M, Sheikholeslami A (2024) Graph clustering with Boltzmann machines. Discret Appl Math 343:208-223. \url{https://doi.org/10.1016/j.dam.2023.10.012}

\bibitem[Spielman et al. 1996]{no17} Spielman D A, Teng S H (1996) Spectral partitioning works: Planar graphs and finite element meshes. In proceedings of 37th conference on foundations of computer science IEEE 1996:96-105. \url{https://doi.org/10.1109/SFCS.1996.548468}

\bibitem[Von Luxburg 2007]{no18} Von Luxburg U (2007) A tutorial on spectral clustering. Stat Comput 17:395-416. \url{https://doi.org/10.1007/s11222-007-9033-z}

\bibitem[Van Dongen 2000]{no19} Van Dongen S (2000) Graph Clustering by Flow Simulation (Ph.D. thesis). Faculteit Wiskunde en Informatica, Universiteit Utrecht.

\bibitem[Hofman et al. 2008]{no20} Hofman J M, Wiggins C H (2008) Bayesian Approach to Network Modularity. Phys Rev Lett 100(25):258701. \url{https://doi.org/10.1103/PhysRevLett.100.258701}

\bibitem[Mackay 2003]{sta-model} Mackay D J C (2003) Information Theory, Inference, and Learning Algorithms. Cambridge University Press, Cambridge, UK. ISBN 0-521-64298-1

\bibitem[Jeub et al. 2015]{no21} Jeub L G S, Balachandran P, Porter M A et al (2015) Think locally, act locally: Detection of small, medium-sized, and large communities in large networks. Phys Rev E 91(1):012821. \url{https://doi.org/10.1103/PhysRevE.91.012821}

\bibitem[Ronhovde et al. 2010]{no22-1} Ronhovde P, Nussinov Z (2010) Local resolution-limit-free potts model for community detection. Phys Rev E 81(4):046114. \url{https://doi.org/10.1103/PhysRevE.81.046114}

\bibitem[Traag et al. 2011]{no22} Traag V A, Van Dooren P, Nesterov Y (2011) Narrow scope for resolution-limit-free community detection. Phys Rev E 84(1):016114. \url{https://doi.org/10.1103/PhysRevE.84.016114}

\bibitem[Boccaletti et al. 2007]{no23} Boccaletti S, Ivanchenko M, Latora V (2007) Detecting complex network modularity by dynamical clustering. Phys Rev E 75(4):045102. \url{https://doi.org/10.1103/PhysRevE.75.045102}

\bibitem[Newman et al. 2004]{no24} Newman M E J, Girvan M (2004) Finding and evaluating community structure in networks. Phys Rev E 69(2):026113. \url{https://doi.org/10.1103/PhysRevE.69.026113}

\bibitem[Fortunato et al. 2004]{no25} Fortunato S, Latora V, Marchiori M (2004) Method to find community structures based on information centrality. Phys Rev E 70(5):056104. \url{https://doi.org/10.1103/PhysRevE.70.056104}

\bibitem[Hagberg et al. 2008]{no26} Hagberg A, Swart P J, Schult D A (2008) Exploring network structure, dynamics, and function using NetworkX. In: Proceedings of the 7th python in science conference, pp 11-15.
\url{https://www.osti.gov/biblio/960616}

\bibitem[Condon et al. 2001]{PPM} Condon A, Karp R M (2001) Algorithms for graph partitioning on the planted partition model. Random Struct Algorithms 18(2):116-140. \url{https://doi.org/10.1002/1098-2418(200103)18:2<116::AID-RSA1001>3.0.CO;2-2}

\bibitem[Fan et al. 2010]{no27} Fan N, Pardalos P M (2010) Linear and quadratic programming approaches for the general graph partitioning problem. J Global Optim 48(1):57-71. \url{https://doi.org/10.1007/s10898-009-9520-1}

\bibitem[Fortunato 2010]{s-1} Fortunato S (2010) Community detection in graphs. Phys Rep 486(3-5):75-174. \url{https://doi.org/10.1016/j.physrep.2009.11.002}

\bibitem[Fortunato et al. 2016]{s-3} Fortunato S, Hric D (2016) Community detection in networks: A user guide. Phys Rep 659:1-44. \url{https://doi.org/10.1016/j.physrep.2016.09.002}

\bibitem[Blondel et al. 2008]{no28} Blondel V D, Guillaume J L, Lambiotte R et al (2008) Fast unfolding of communities in large networks. J Stat Mech-Theory Exp 2008(10):P10008. \url{https://doi.org/10.1088/1742-5468/ad6139}

\bibitem[Fortunato et al. 2007]{no29} Fortunato S, Barthelemy (2007) Resolution limit in community detection. In proceedings of the national academy of sciences, 104(1):36-41. \url{https://doi.org/10.1073/pnas.0605965104}

\bibitem[Good et al. 2010]{no30} Good B H, De Montjoye Y A, Clauset A (2010) Performance of modularity maximization in practical contexts. Phys Rev E 81(4):046106. \url{https://doi.org/10.1103/PhysRevE.81.046106}

\bibitem[Holland et al. 1983]{SBM} Holland P W, Laskey K B, Leinhardt S (1983) Stochastic blockmodels: First steps. Soc Networks 5(2):109-137. \url{https://doi.org/10.1016/0378-8733(83)90021-7}

\bibitem[Wolsey 2020]{Integer-programming} Wolsey L A (2020) Integer programming. John Wiley \& Sons. ISBN 9781119606536


\bibitem[Fan et al. 2012]{no33} Fan N, Zheng Q P, Pardalos P M (2012) Robust optimization of graph partitioning involving interval uncertainty. Theoret Comput Sci 447:53-61. \url{https://doi.org/10.1016/j.tcs.2011.10.015}

\bibitem[Miasnikof et al. 2020a]{no34} Miasnikof P, Pitsoulis L, Bonner A J et al (2020) Graph clustering via intra-cluster density maximization. In: Network Algorithms, Data Mining, and Applications: NET, Moscow, Russia, May 2018 8. Springer International Publishing, pp 37-48. \url{https://doi.org/10.1007/978-3-030-37157-9_3}

\bibitem[Ponomarenko et al. 2021]{no35} Ponomarenko A, Pitsoulis L, Shamshetdinov M (2021) Overlapping community detection in networks based on link partitioning and partitioning around medoids. PLoS One 16(8):e0255717. \url{https://doi.org/10.1371/journal.pone.0255717}

\bibitem[Zhao et al. 2021]{no36} Zhao P F, Li Q N, Chen W K, et al (2021) An efficient quadratic programming relaxation based algorithm for large-scale MIMO detection. SIAM J Optim 31(2):1519-1545. \url{https://doi.org/10.1137/20M1346912}

\bibitem[Cui et al. 2018]{no37} Cui C, Li Q, Qi L, et al (2018) A quadratic penalty method for hypergraph matching. J Glob Optim 70:237-259. \url{https://doi.org/10.1007/s10898-017-0583-0}

\bibitem[Miasnikof et al. 2020b]{cq-1} Miasnikof P, Shestopaloff A Y, Bonner A J et al (2020) A density-based statistical analysis of graph clustering algorithm performance. J Complex Netw 8(3):cnaa012. \url{https://doi.org/10.1093/comnet/cnaa012}

\bibitem[Miasnikof et al. 2018]{cq-2} Miasnikof P, Shestopaloff A Y, Bonner A J et al (2018) A statistical performance analysis of graph clustering algorithms. In: Algorithms and Models for the Web Graph: 15th International Workshop, WAW 2018, Moscow, Russia, May 17-18, 2018, Proceedings 15. Springer International Publishing, pp 170-184. \url{https://doi.org/10.1007/978-3-319-92871-5_11}

\bibitem[Newman 2004b]{cq-3} Newman M E J (2004) Analysis of weighted networks. Phys Rev E 70(5):056131. \url{https://doi.org/10.1103/PhysRevE.70.056131}

\bibitem[Burt 1976]{Burt-distance} Burt R S (1976) Positions in networks. Social forces 55(1):93-122. \url{https://doi.org/10.1093/sf/55.1.93}
\bibitem[Jaccard 1901]{jaccard-1} Jaccard P (1901) \'{e}tude de la distribution florale dans une portion des Alpes et du Jura. Bull Soc Vaudoise Sci Nat 37:547-579. 

\bibitem[Camby et al. 2017]{jaccard-2} Camby \'{E}, Caporossi G (2017) The extended Jaccard distance in complex networks. GERAD, \'{E}cole des hautes \'{e}tudes commerciales.

\bibitem[2022]{jaccard-3} Miasnikof P, Shestopaloff A Y, Pitsoulis L et al (2022) An empirical comparison of connectivity-based distances on a graph and their computational scalability. J Complex Netw 10(1):cnac003. \url{https://doi.org/10.1093/comnet/cnac003}

\bibitem[Miasnikof et al. 2021]{jaccard-4} Miasnikof P, Shestopaloff A Y, Pitsoulis L et al (2021) Distances on a graph. In: Complex Networks \& Their Applications IX: Volume 1, Proceedings of the Ninth International Conference on Complex Networks and Their Applications COMPLEX NETWORKS 2020. Springer International Publishing, pp 189-199. \url{https://doi.org/10.1007/978-3-030-65347-7_16}

\bibitem[Garcia-Mendez et al. 2020]{jaccard-5} Garcia-Mendez S, Fernandez-Gavilanes M, Juncal-Martinez J et al (2020) Identifying banking transaction descriptions via support vector machine short-text classification based on a specialized labelled corpus. IEEE Access 8(2020):61642-61655. \url{https://doi.org/10.1109/ACCESS.2020.2983584}

\bibitem[Wang et al. 2020]{jaccard-6} Wang Z B, Cui J, Zhu Y (2020) Plant recognition based on Jaccard distance and BOW. Multimedia Syst 26(5):495-508. \url{https://doi.org/10.1007/s00530-020-00657-6}

\bibitem[Ochiai 1957]{O-O-distance} Ochiai A (1957) Zoogeographical studies on the soleoid fishes found in japan and its neighbouring regions-i. NIPPON SUISAN GAKKAISHI 22(9):522-525. \url{https://cir.nii.ac.jp/crid/1571980075647016960}

\bibitem[Nocedal et al. 2006]{nno5} Nocedal J, Wright S (2006) Numerical Optimization. Springer, Berlin. ISBN-13 978-0387-30303-1

\bibitem[Sun et al. 2006]{nno6} Sun W, Yuan Y X (2006) Optimization theory and methods: nonlinear programming. Springer Science \& Business Media, New York. ISBN-10 0-387-24975-3

\bibitem[Birgin et al. 2000]{NSPG} Birgin E G, Mart\'{i}nez J M and Raydan M (2000) Nonmonotone spectral projected gradient methods on convex sets. SIAM J Optim, 10(2000):1196-1121.
\url{https://doi.org/10.1137/S1052623497330963}


\bibitem[Barzilai et al. 1988]{BB3} Barzilai J, Borwein J M (1988) Two point step size gradient methods. IMA J Numer Anal, 8 (1988):141-148.
\url{https://doi.org/10.1093/imanum/8.1.141}

\bibitem[Jeub et al. 2011-2019]{louvain} Jeub L, Bazzi M, Jutla I, Mucha P (2011-2019)
A generalized Louvain method for community detection implemented in MATLAB. \url{https://github.com/GenLouvain/GenLouvain}

\bibitem[Gurobi 2023]{Gurobi} Gurobi Optimization, LLC, Gurobi Optimizer Reference Manual. 2023. [Online]. Available: \url{https://www.gurobi.com}

\end{thebibliography}

\end{sloppypar}
\end{document}